\documentclass[a4paper, 12pt]{article}
\usepackage{amsmath}

\newtheorem{thm}{Theorem}
\newtheorem{prop}{Proposition}
\newtheorem{lemma}{Lemma}
\newtheorem{cor}{Corollary}
\newtheorem{defn}{Definition}

\newcommand{\Haus}{\mathcal{H}}

\newcommand{\deltat}{\tilde{\delta}}

\newcommand{\Rm}{\mathbf{R}^m}
\newcommand{\Rmm}{\mathbf{R}^{2m}}
\newcommand{\Rnn}{\mathbf{R}^{2n}}
\newcommand{\coneprod}{\times\!\!\!\!\!\times}
\newcommand{\be}{\begin{equation}}
\newcommand{\ee}{\end{equation}}
\newcommand{\bea}{\begin{eqnarray}}
\newcommand{\eea}{\end{eqnarray}}
\DeclareMathOperator{\dive}{div}
\DeclareMathOperator{\spt}{spt}
\DeclareMathOperator{\Lip}{Lip}
\DeclareMathOperator{\dist}{dist}
\newcommand\res{\mathop{\hbox{\vrule height 6pt width .5pt depth 0pt
\vrule height .5pt width 6pt depth 0pt}}\nolimits}
\newcommand\ins{\mathop{\hbox{\vrule height .5pt width 6pt depth 0pt
\vrule height 6pt width .5pt depth 0pt}}\nolimits}
\newcommand\avint{\hbox{\hbox{$\displaystyle \int$}\hbox{\kern-.9em{$-$}}}}

\begin{document}

\title{Uniqueness of tangent cones for calibrated $2$-cycles}
\author{David Pumberger\footnote{Partly supported by the Swiss National Fund} \\ D-Math\\ ETH Z\"urich\\ CH-8092 Z\"urich\\ Switzerland\\ maupu@math.ethz.ch \and Tristan Rivi\`ere\\ D-Math\\ ETH Z\"urich\\ CH-8092 Z\"urich\\ Switzerland\\ riviere@math.ethz.ch}
\date{}
\maketitle
\begin{abstract}
In this paper we prove that the tangent cones to calibrated $2$-cycles are unique. Furthermore, using this result we prove a rate of convergence for the mass of the blow-up of a calibrated integral $2$-cycle $C$ towards the limiting density: there exist constants $C_1>0$, $\gamma >0$ such that
$$
\frac{M(C\res B_r(x_0))}{r^2}-\Theta (\|C\|, x_0)\leq C_1 r^\gamma\ .
$$
We also obtain such a rate for $J$-holomorphic maps between almost complex manifolds and deduce that their tangent maps are unique.
\end{abstract}


\section{Introduction}

Let $M$ be a smooth compact $m$-dimensional manifold and let $\Omega _0^k(M)$ denote the smooth compactly supported $k$-forms on $M$. Then a $k$-dimensional current $C$ in $M$ is a distribution on the compactly supported $k$-forms on $M$. The boundary of such a $k$-current is the $(k-1)$-current defined by $\partial C(\omega):=C(d\omega)$, where $\omega\in \Omega _0^{k-1}(M)$ and we say that a $k$-current $C$ is a $k$-cycle if $\partial C=0$.\\
We can put an arbitrary smooth Riemannian metric $g$ on $M$ and define the comass of a $k$-form $\omega$ to be 
$$
\|\omega \|_*:=\sup _{x\in M}\sup _{e_1,\ldots , e_k\in S_xM}|\langle \omega, e_1\wedge\ldots \wedge e_k\rangle |\ ,
$$
where $S_xM$ means the unit sphere in $T_xM$ with respect to $g_x$. Then the mass $M(C)$ of a $k$-current $C$ is defined by 
$$
M(C):=\sup _{\omega\in \Omega _0^k(M), \|\omega \|_*\leq 1} |C(\omega)|\ .
$$
A $k$-cycle $C$ is called a normal cycle if $C$ satisfies $M(C)<+\infty$. Furthermore, we call a current an integer multiplicity rectifiable $k$-current, if $C$ is a normal $k$-current and there are $k$-Hausdorff measurable subsets $\mathcal{N}_j$ of oriented $k$-dimensional $C^1$-submanifolds $N_j$ with $\mathcal{N}_i\cap \mathcal{N}_j=\emptyset$, $i\neq j$, and a multiplicity $\Theta:\bigcup _{j=1}^\infty \mathcal{N}_j\to\mathbf{Z}$ such that
$$
\langle C; \psi \rangle =\sum _{j=1}^\infty \int \psi\, \Theta\ d\mathcal{H}^k\res\mathcal{N}_j\ .
$$
In this case the mass of $C$ is given by $M(C)=\sum _{j=1}^\infty \int |\Theta |\ d\mathcal{H}^k\res\mathcal{N}_j$. In this paper we will only work with $k$-cycles and abbreviate integer multiplicity rectifiable $k$-cycles by calling them integral $k$-cycles (usually integral currents are currents for which both $C$ and $\partial C$ are integer multiplicity rectifiable currents).\\

We also need to define the notion of a smooth calibration on $M$:
\begin{defn}
Let $\omega$ be a smooth $k$-form on $(M^m,g)$. Then we say that a $k$-current is calibrated by $\omega$ if for all open subsets $U\subset M$ we have 
$$
\langle C\res U; \omega\rangle =M(C\res U)\ ,
$$
where $C\res U$ means the restriction of $C$ to $U$. 
\end{defn}
If $\omega$ has comass equal to $1$ and satisfies $d\omega =0$, it is easy to see, that a current calibrated by $\omega$ is homologically mass-minimizing (see Harvey and Lawson \cite{hala82} for details). In this paper we will always assume our calibrations to have comass equal to $1$ but for reasons which will become apparent later, we do not assume the calibration to be closed. Thus calibrated currents will not always be (homologically) mass-minimizing.\\
We now want to give some well-known and important examples of calibrated currents. In complex geometry (complex) $p$-dimensional algebraic sub-varieties of $\mathbf{CP}^n$ can be seen as $2p$-dimensional integral currents calibrated by $\omega _{\mathbf{CP}^n}^p$, where $\omega _{\mathbf{CP}^n}$ is the standard K\"ahler form on $\mathbf{CP}^n$. More generally, when $(M^{2m}, J,\omega)$ is an almost K\"ahler manifold where $\omega$ is the closed K\"ahler form compatible with $J$ (this means that $\omega (\cdot, J\cdot)$ is a Riemannian metric on $M$), integral $p-p$-currents are calibrated by $\omega ^p$ --- here a $p-p$-cycle is an integral cycle where the $C^1$-manifolds $N_j$ satisfy $J_x(T_xN_j)=T_xN_j$ for any $x\in N_j$. Note that in this case $p-p$-currents are a subclass of the area-minimizing $2p$-currents.\\
On a general compact almost complex manifold $(M^{2m}, J)$ one can still define $p-p$-currents as above. Notice that the case of $1-1$-cycles is particularly interesting, since they arise as perturbations of $J_0$-holomorphic graphs and are generic from the existence point of view --- see the introduction of \cite{riti03}. If one can find a symplectic form which is compatible with $J$, then we are back in the almost K\"ahler case and $1-1$-cycles are area-minimizing. It turns out that locally one can indeed find a $2$-form $\omega$ which is compatible with $J$ but in general this $\omega$ will not be closed (hence not symplectic). In real dimensions $\leq 4$ the form $\omega$ can be constructed to be closed --- see the appendix of \cite{riti01} --- but for higher dimensions R. Bryant \cite{br82} constructed an almost complex structure on $S^6$ which does not admit any compatible $\omega$ even locally. Hence in such a case $1-1$-cycles are still calibrated by $\omega$ but no longer area-minimizing. However, $1-1$-cycles are of geometric interest even in the absence of a symplectic form when results for area-minimizing currents no longer apply. This is one of the reasons, why we do not want to assume the calibration to be closed in this paper.\\
Another important example of calibrated currents is connected to Special Lagrangian currents in Calabi-Yau $p$-folds. A Calabi-Yau $p$-fold is a $2p$-dimensional K\"ahler manifold $(M^{2p}, J, \omega)$ for which the holonomy is $SU(p)$. As one can see in the lecture notes by D. Joyce \cite{jo01} , for such a manifold there exists a $(p,0)$-holomorphic form $\Omega$ (called the holomorphic volume form) satisfying
$$
\frac{\omega ^p}{p!}=(-1)^\frac{p(p1)}{2} \bigg(\frac{i}{2}\bigg)^p \Omega\wedge\bar{\Omega}\ .
$$
Then a $p$-current $C$ is called Special Lagrangian if it is calibrated by $\text{Re}\, \Omega$ --- the real part of $\Omega$. Tangent cones $C_\infty$ to such currents (this will be explained below) are calibrated by $\text{Re}(dz_1\wedge\ldots \wedge dz_p)$ (where we use usual normal coordinates) and are of the form $C_\infty =0\coneprod T$ (see \cite{fe69} for the notation), where $T$ is a Special Legendrian $p-2$-current in $S^{2p-2}$. These $p-2$-currents are calibrated by $\text{Re}\, (\sum z_i dz_{i+1}\wedge dz_{i-1})$ and it is important to note that again this calibration is not closed.\\

Our goal is to study the regularity of such calibrated cycles so it is natural to analyze their blow-up around points first, i.e. to prove whether such a blow-up exists or not and when it exists whether it is unique or not. The blow-up analysis of $C$ around a point $x_0\in M$ is done as follows: consider a dilation of $C$ around $x_0$ which in normal coordinates near $x_0$ is given by the push-forward of the current $C$ under the map $\frac{x-x_0}{r}$ --- here we mean that $\frac{x-x_0}{r}_*C(\psi):=C(\frac{x-x_0}{r}^*\psi)$. To analyze the behavior of these dilations as $r\to 0$ we need the monotonicity formula which says that $r^{-k} M(C\res B_r(x_0))$ is increasing and the compactness theorem by Federer-Fleming (see \cite{fe69} $4.2.17$). The monotonicity formula gives that $M(\frac{x-x_0}{r}_*C\res B_1(0))$ is bounded independent of $r$ since $M(\frac{x-x_0}{r}_*C\res B_1(0))\cong r^{-k} M(C\res B_r(x_0))$. Using this and the fact that $C$ is a cycle we apply the compactness theorem to deduce that there exists $r_n\to 0$ and a normal current $C_\infty$ such that weakly
$$
\frac{x-x_0}{r_n}_*C\res B_1\rightharpoonup C_\infty\res B_1\ .
$$
It turns out that $C_\infty$ is a cone --- called a tangent cone to $C$ at $x_0$ --- which is calibrated by $\omega _{x_0}$ (see section \ref{preliminaries} of this paper). The main questions related to such a construction are whether the blow-up limit is unique or not and whether the dilated currents converge to the limiting object at a certain rate.\\

These questions are directly related to rate of convergence results for elliptic systems obtained by C. Morrey (see \cite{mo66} or \cite{giaq83}). For example, if one considers stationary harmonic maps $u:B^m\to N^n\subset\mathbf{R}^N$ from the unit ball in $\mathbf{R}^m$ into a submanifold $N$ of $\mathbf{R}^N$ --- for a detailed discussion we refer to the books by F. H\'elein \cite{he02} and L. Simon \cite{si96} ---, one obtains that $\lim_{r\to 0}\frac{1}{r^{m-2}}\int _{B_r(x_0)}|\nabla u|^2=:\Theta _u(x_0)$ exists for all $x_0\in B^m$. This implies that the dilated maps $u_{x_0,r}(x):=u(x_0+rx)$ are uniformly bounded in $W^{1,2}$ and hence that there exists a sequence $r_n\to 0$ and a $W^{1,2}$-map $u_{x_0,\infty}$ --- a tangent map to $u$ at $x_0$ --- such that $u_{x_0,r_n}\rightharpoonup u_{x_0,\infty}$. In case $\Theta _u(x_0)=0$, C. Evans \cite{ev91} and F. Bethuel \cite{be93} showed that there exist constants $C_1>0$, $\gamma >0$ such that 
$$
\frac{1}{r^{m-2}}\int _{B_r(x_0)} |\nabla u|^2 \leq C_1 r^\gamma\ .
$$
Using results of Morrey, this rate implies that in a small neighborhood $U$ of $x_0$ we have that $u\in C^{0,\frac{\gamma}{2}}(U,N)$. Then the $u_{x_0,r}$ converge to a unique tangent map (independent of the sequence $r_n$) $u_{x_0,\infty}$ which is constant. \\ However, it would also interesting to know, if for $\Theta _u(x_0)>0$ there are constants $C_1>0$, $\gamma >0$ such that the following more general estimate holds:
 $$
\frac{1}{r^{m-2}}\int _{B_r(x_0)} |\nabla u|^2 -\Theta _u(x_0) \leq C_1 r^\gamma\ . 
 $$
In other words, one would like to have a rate of convergence for $\frac{1}{r^{m-2}}\int _{B_r(x_0)} |\nabla u|^2$ towards the limiting density $\Theta _u(x_0)$. Note that again such a rate of convergence is very strong information as it implies the weak convergence of the $u_{x_0,r}$ to a unique blow-up limit $u_{x_0,\infty}$ as shown in proposition \ref{mapsunique}.\\
Unfortunately, Gulliver and White showed in \cite{guwh89} that such a rate does not always hold for stationary harmonic maps. However, in the first section of this paper we show, that such a rate in fact exists for locally approximable $J$-holomorphic maps. For a compact almost complex manifold $(M, J_M)$ and a tamed compact symplectic manifold $(N,J_N, \omega _N)$ a $W^{1,2}$-map $u:M\to N$ is called a locally approximable $J$-holomorphic map, if for a.e. $x\in M$ and all $X\in T_xM$ we have $du (J_M(X))=J_N(du(X))$ and $u$ locally is in the strong closure of the $C^\infty$-maps in $W^{1,2}(M,N)$ (note that in special cases $J$-holomorphic maps are stationary harmonic --- see \cite{riti01} or \cite{wa03}). For such a map we can prove the following theorem:

\begin{thm}\label{J-maps}
Let $u:M\rightarrow N$ be a locally approximable $J$-holomorphic map from an almost complex manifold $(M, J)$ into a tamed compact symplectic manifold $(N, J_N, \omega _N)$, such that $u\in W^{1,2}(M,N)$. Given any $x_0\in M$ there exist $C_1$, $C_2>0$ and $\gamma\in (0,1]$ independent of $u$ such that 
$$\frac{1}{r^{2n-2}}\int _{B_r(x_0)}|\nabla u|^2 -\Theta _u(x_0)\leq C_1\big( 1+C_2\| \nabla u\| _{L^2(B_1)}^2 \big) r ^\gamma\ ,$$ 
where $\Theta _u(x_0):=\lim _{\rho\rightarrow 0}\frac{1}{\rho ^{2n-2}}\int _{B_\rho(x_0)} |\nabla u|^2$.
\end{thm}

Now we go back to the case of area-minimizing integral $k$-cycles. There B. White showed that a rate of convergence for $r^{-k}M(C\res B_r(x_0)$ to the limiting density $\Theta (\|C\|, x_0)$ implies the uniqueness of the tangent cone to $C$ at $x_0$ --- see theorem 3 in \cite{wh83} ---, i.e. if there exist constants $C_1>0$ and $\gamma>0$ such that 
$$
\frac{M(C\res B_r(x_0))}{r^k}-\Theta (\|C\|, x_0)\leq C_1 r^\gamma\ ,
$$
then the blow-up limit is independent of the subsequence chosen. Using an idea of E. Reifenberg he deduced such a rate from a comparison argument in the case of area-minimizing integral $2$-cycles. In section \ref{uniqueness} of this paper we will show that the tangent cones to normal $2$-cycles calibrated by a not necessarily closed form are unique:

\begin{thm}\label{uniqueness}
Let $C$ be a normal $2$-cycle in $\mathbf{R}^m$ which is calibrated by a $C^2$-form $\omega$ (not necessarily closed). Then for any $x_0\in\mathbf{R}^m$ there exists a unique tangent cone $C_\infty$ to $C$ at $x_0$. 
\end{thm}

The proof we give only depends on geometric observations and the monotonicity formula, so we do not use the approach through a comparison argument. The monotonicity we use will be proved in \ref{calmon}. It is similar to the one for area-minimizing integral $k$-cycles obtained by H. Federer in chapter 5 of \cite{fe69} (5.4.3 (2)) and was proven for calibrated currents (in case the calibration is a closed form) by R. Harvey and B. Lawson \cite{hala82}. In fact the monotonicity formula we obtain is only an almost monotonicity formula but turns out to be good enough to study blow-ups.\\ 
We then use theorem \ref{uniqueness} to deduce a rate of convergence for integral $2$-cycles:

\begin{thm}\label{rate}
If, in addition to the hypotheses of theorem \ref{uniqueness}, $C$ is also integer rectifiable, then for any $x_0\in\mathbf{R}^m$ there exist $r_0>0$, $C_1>0$ and $\gamma >0$ such that for $\rho\in (0,r_0)$ we have 
$$
\frac{M(C\res B_\rho(x_0))}{\rho ^2}-\Theta (\|C\|, x_0) \leq C_1\rho^\gamma\ .
$$
\end{thm}

As we already mentioned above, uniqueness of tangent cones has been proven in several related situations --- for area-minimizing integral $2$-cycles we refer the reader to the paper by B. White \cite{wh83}. The problem of uniqueness was also studied by W. Allard and F. Almgren in \cite{alal76} and \cite{alal81}, L. Simon in \cite{si83} and J. Taylor in \cite{ta73} and \cite{ta76}.\\

Our paper is organized as follows: in section \ref{maps} we prove a rate of convergence for $J$-holomorphic maps and deduce the uniqueness of their tangent maps. In section \ref{currents} we prove theorems \ref{uniqueness} and \ref{rate} first for constant calibrations and then for general calibrations through a perturbation argument.


\section{Proof of theorem \ref{J-maps}}\label{maps}

Since the argument is local, we can assume that the domain manifold is in fact the open unit ball  $B^{2n}\subset\mathbf{R}^{2n}$ and that $x_0=0\in B^{2n}$. Also assume that the target is isometrically embedded in $\mathbf{R}^{2N}$ for some $N$ large enough (see \cite{pu04} for details). To make the presentation clearer we split the proof of theorem \ref{J-maps} into two cases. First we prove the theorem in the case where the almost complex structure is the standard complex structure $J_0$ of $\mathbf{R}^{2n}$ to illustrate the ideas. The remaining case will be handled by a perturbation argument.

\subsection{The case of the standard complex structure}

In this case we know that $J$-holomorphic maps are stationary harmonic (in fact energy minimising in their homotopy class --- see for instance \cite{riti01}). Therefore we have the usual monotonicity formula: 
$$
\frac{d}{dr}\bigg(\frac{1}{r^{2n-2}}\int _{B_r(x_0)} |\nabla u|^2\bigg)=2\int _{\partial B_r(x_0)}R^{2-2n}\Big|\frac{\partial u}{\partial R}\Big|^2\ .
$$ 
Integrating this formula from $0$ to $r$ we obtain that 
$$
2\int _0^r \frac{d\rho}{\rho^{2n-2}}\int _{\partial B_\rho(x_0)} \Big|\frac{\partial u}{\partial R}\Big|^2 =\frac{1}{r^{2n-2}}\int _{B_r (x_0)} |\nabla u|^2 -\Theta _u(x_0)<+\infty\ 
$$ 
i.e. that $\int _{B_r}|x-x_0|^{2-2n}|\frac{\partial u}{\partial R}|^2<+\infty$. Next we want to show that given almost any $J$-holomorphic plane in $\mathbf{R}^{2n}$ the restriction of $u$ to such a plane is a $J$-holomorphic $W^{1,2}$-map. To see this first note that given any unit vector $X$ it lies in the plane spanned by $\{X,JX\}$ and that 
$$
\int _{B_r(x_0)}|x-x_0|^{2-2n}|u\cdot X|^2 \leq \int _{B_r (x_0)}|x-x_0|^{2-2n}\Big|\frac{\partial u}{\partial R}\Big|^2 <+\infty 
$$
by the monotonicity formula. For the slicing let $f$ be the map parametrising the $J$-holomorphic curves through $x_0$, i.e. $f:B^{2n}(x_0)\rightarrow \mathbf{CP}^{n-1}$. Then this map satisfies $|\nabla f|=O(\frac{1}{r})$, so that there are constants $C_1, C_2>0$ with $C_1-C_2|x-x_0| \leq |x-x_0|^{2n-2} |J_{2n-2}f|\leq C_1+C_2|x-x_0|$ (here $J_{2n-2}f$ denotes the $(2n-2)$-Jacobian of $f$). We use this fact when we now apply the co-area formula: 
$$
\int _{\mathbf{CP}^{n-1}}dp\int _{f^{-1}(p)}\frac{|x-x_0|^{2-2n}|u\cdot X|^2}{|J_{2n-2}f|}=\int _{B^{2n}(x_0)}|x-x_0|^{2-2n}|u\cdot X|^2 <+\infty\ .
$$ 
By Fubini's theorem this implies that for almost every $p\in\mathbf{CP}^{n-1}$, $\int _{f^{-1}(p)}|u\cdot X|^2 <+\infty$. Since the map $u$ and $f^{-1}(p)$ are $J$-holomorphic we conclude that for such $p$ we also have $\int _{f^{-1}(p)}|u\cdot JX|^2 <+\infty$ and hence that for almost every $p\in\mathbf{CP}^{n-1}$, $\int _{f^{-1}(p)}|\nabla u|^2 <+\infty$. Therefore $u|_{f^{-1}(p)}$ is a $J$-holomorphic $W^{1,2}$-map from $(f^{-1}(p),J_0)$ to $(N,J_N)$. From the regularity theorems for $J$-holomorphic maps from Riemann surfaces we conclude that the restricted map is smooth and that there exists $\delta >0$ such that 
$$
\int _{f^{-1}(p)\cap B_r(x_0)}|\nabla u|^2 <\delta \int _{f^{-1}(p)\cap B_{2r}(x_0)}|\nabla u|^2\ .
$$
Integrating this estimate over $p\in\mathbf{CP}^{n-1}$ yields 
$$
\int _{\mathbf{CP}^{n-1}}\int _{f^{-1}(p)\cap B_r(x_0)} |u\cdot X|^2 <\delta \int _{\mathbf{CP}^{n-1}}\int 
_{f^{-1}(p)\cap B_{2r}(x_0)} |u\cdot X|^2
$$ 
since $|\nabla u|^2=|u\cdot X|^2+|u\cdot JX|^2=2|u\cdot X|^2$. Using the co-area formula again this gives
\begin{eqnarray} 
\nonumber \int _{\mathbf{CP}^{n-1}}\int _{f^{-1}(p)\cap B_r(x_0)} |u\cdot X|^2 & = & \int _{B_r(x_0)} |u\cdot X|^2 |J_{2n-2}f| \\ 
\nonumber & < & \delta \int _{B_r(x_0)} |u\cdot X|^2 |J_{2n-2}f|\ .
\end{eqnarray}
Next we take $r>0$ sufficiently small so that $\frac{C_1+C_2r}{C_1-C_2r}$ is close enough to $1$ for  $\frac{C_1+C_2r}{C_1-C_2r}\delta <1$. Then for some $\deltat \in (0,1)$ the above estimate becomes 
$$
\int _{B_r(x_0)}|x-x_0|^{2-2n}|u\cdot X|^2 < \deltat \int _{B_{2r}(x_0)}|x-x_0|^{2-2n} |u\cdot X|^2\ .
$$ 
Therefore we get that there exist $\gamma \in (0,1]$ and $C>0$ such that 
$$
\int _{B_r(x_0)}|x-x_0|^{2-2n} |u\cdot X|^2 \leq Cr^\gamma\ .
$$ 
Furthermore, as above, we get that 
$$
\int _{B_r(x_0)}|x-x_0|^{2-2n}\Big|\frac{\partial u}{\partial R}\Big|^2\leq \int _{B_r(x_0)}|x-x_0|^{2-2n} |u\cdot X|^2 \leq Cr^\gamma 
$$ 
which together with the monotonicity formula yields 
\begin{eqnarray} 
\nonumber \frac{1}{2}\int _0^r d\rho \frac{d}{d\rho}\bigg(\rho^{2-2n}\int _{B_\rho (x_0)} |\nabla u|^2\bigg) & = & \int _0^r d\rho \int _{\partial B_\rho (x_0)}R^{2-2n}\Big|\frac{\partial u}{\partial R}\Big|^2\\ 
\nonumber & = & \int _{B_r(x_0)}|x-x_0|^{2-2n}\Big|\frac{\partial u}{\partial R}\Big|^2\\ 
\nonumber & \leq & Cr^\gamma \ , 
\end{eqnarray}
i.e. 
$$
\frac{1}{r^{2n-2}}\int _{B_r(x_0)} |\nabla u|^2 -\Theta _u(x_0)\leq Cr^\gamma\ ,
$$ 
which is the desired estimate in the easy case.

\subsection{Perturbing the proof}

In this subsection we show that the above proof can still be carried out when the domain of $u$ is the unit ball $B^{2n}\subset\mathbf{R}^{2n}$ together with an arbitrary complex structure $J$ (the proof for an arbitrary metric is similar). Without loss of generality (through a change of coordinates) we can always assume that this complex structure is the standard complex structure at the origin, i.e. that $J(0)=J_0$, and that $J$ satisfies $\| J(r)-J_0\|_{C^2}=O(r)$ on the whole ball. We will show that these assumption suffice to carry out the above proof. To do this we first prove that a locally approximable $J$-holomorphic map satisfies an equation similar to the equation for stationary harmonic maps (see for example L. Simon \cite{si96}). We then proceed to deduce an almost monotonicity formula which will give us that the density exists at any point. Such a monotonicity formula was first proved by P. DeBartholomeis and G. Tian \cite{dbti96}. However, they did not show the precise error term which is crucial  in the proof below. We proceed by showing that all the arguments in the $J_0$-case can still be carried out introducing only a small extra error.

\subsubsection{A second order elliptic equation for $u$}
We begin by showing that a locally approximable $J$-holomorphic map is almost a critical point of the  Dirichlet energy $E(u):=\frac{1}{2}\int _B |\nabla u|^2$ for perturbations in the domain. Given an arbitrary smooth $1$-parameter family $F_t$ of diffeomorphisms of $B^{2n}$ we want to compute $\frac{d}{dt} E(u\circ F_t)\bigg|_{t=0}$. First note that one can easily verify the following alternative expression for the energy:
$$
E(u)=\frac{1}{4}\int _B |\nabla u+J_N\circ\nabla u\circ J_0|^2+\frac{1}{2}\int _B\langle J_N\nabla u,\nabla u\circ J_0\rangle\ .
$$
To simplify the notation we write $u_t:=u\circ F_t$. Note that for any map $u_t$ in the second term above we can write $\langle J_N\circ \nabla u_t,\nabla u_t\circ J_0\rangle =\langle \omega _0, u_t^*\omega _N\rangle$ where $\omega _0$, $\omega _N$ are the symplectic structures compatible with the metrics on domain and target respectively. For smooth perturbations of the domain we claim that the second term in the above expression satisfies
$$
\frac{d}{dt} \int _{B^{2n}} \langle \omega _0, u_t^*\omega _N\rangle \bigg|_{t=0}=0\ .
$$
To see this note that on $B^{2n}$ there exists a smooth $1$-form $\phi$ such that $d\phi =\omega _0$ and that $u_t^*\omega _N=F_t^*u^*\omega _N$, whence $d(u_t^*\omega _N)=F_t^*d(u^*\omega _N)=0$ as $u$ is locally approximable. Using Stoke's theorem and the fact that $u_t$ has compact support in $B^{2n}$ we get  
\begin{eqnarray}
\nonumber \int _{B^{2n}} \langle \omega _0, u_t^*\omega _N\rangle & = & \int _{B_{2n}}\langle d\phi, u_t^*\omega _N\rangle \\ 
\nonumber & = & C\int _{B^{2n}} d(\phi\wedge\omega _0^{2n-2})\wedge u_t^*\omega _N\\ 
\nonumber & = & C\int _{B^{2n}} d(\phi\wedge\omega _0^{2n-2}\wedge u_t^*\omega _N)\\ 
\nonumber & = & 0\ .
\end{eqnarray}
Thus it remains to compute the derivative for the first term. We will show that since $J$ is close to $J_0$ in $B^{2n}$, $u$ is close to being $J_0$-holomorphic. Note that 
$$
\nabla u+J_N\circ \nabla u \circ J_0=\nabla u+J_N\circ \nabla u\circ J+J_N\circ \nabla u\circ (J-J_0)
$$
so that for $u_t$ (using the fact that $J_N$ is compatible with $g_N$) this means that 
\begin{eqnarray}
\label{term1} \int _{B^{2n}} |\nabla u_t+J_N\circ \nabla u_t\circ J_0|^2 & = & \int _{B^{2n}} |\nabla u_t+J_N\circ \nabla u_t\circ J|^2\\ 
\label{term2} & + &  2\int _{B^{2n}}\langle \nabla u_t, J_N\circ \nabla u_t\circ (J-J_0)\rangle \\ 
\label{term3} & + & \int _{B^{2n}}\langle \nabla u_t\circ (J_0+J), \nabla u_t \circ (J-J_0)\rangle \ .
\end{eqnarray}
For the first term on the right-hand side (\ref{term1}) note that since $u$ is $J$-holomorphic we have $|\nabla u+J_N\circ \nabla u\circ J_0|^2=0$ pointwise a.e. in $B^{2n}$, i.e. $\int _{B^{2n}} |\nabla u+J_N\circ \nabla u\circ J_0|^2=0$. Also
$$
\int _{B^{2n}} |\nabla u_t+J_N\circ \nabla u_t\circ J|^2\geq \int _{B^{2n}} |\nabla u+J_N\circ \nabla u\circ J|^2=0\ ,
$$
which implies that
$$
\frac{d}{dt}\int _{B^{2n}} |\nabla u_t+J_N\circ \nabla u_t\circ J|^2\bigg|_{t=0} =0\ .
$$
For the remaining terms (\ref{term2}) and (\ref{term3}) we work in local coordinates. Since any $1$-parameter group $F_t$ is generated by a vector field $\xi$ having compact support in $B^{2n}$, we can write $u_t(x)=u(x+t\xi)$. Writing $(J-J_0)(x)=[a_{kl}](x)$ and $J_N(u(x))=[b_{\alpha\beta}](u(x))$ --- throughout the paper we use latin indices to denote coordinates on the domain and greek ones for the target --- we get the following expression for the derivative of the term in (\ref{term2}):
\begin{eqnarray}
\nonumber \frac{d}{dt}\int _{B^{2n}} \langle \nabla u_t, J_N\circ \nabla u_t\circ (J-J_0)\rangle \bigg| _{t=0} & = & \\
\label{term21} & \hspace{-3 cm} = & \hspace{-1,5 cm} \int _{B^{2n}} (-\dive {\xi})\frac{\partial u^\alpha}{\partial x_k}\frac{\partial u^\beta}{\partial x_l}[a_{kl}][b_{\beta\alpha}] \\ 
\label{term22} & \hspace{-3 cm} + & \hspace{-1,5 cm} 2\int _{B^{2n}} \frac{\partial u^\alpha}{\partial x_k}\frac{\partial u^\beta}{\partial x_{l'}}\frac{\partial \xi^{l'}}{\partial x_l}[a_{kl}][b_{\beta\alpha}] \ .
\end{eqnarray}
For the last term (\ref{term3}) we write $(J+J_0)(x)=([a_{kl}](x)+2[\deltat _{kl}])$, where $[\deltat _{kl}]$ denotes the matrix for $J_0$ in the standard coordinates on $B^{2n}$. Then a computation shows that
\begin{eqnarray}
\nonumber \frac{d}{dt} \int _{B^{2n}}\langle \nabla u_t\circ (J_0+J), \nabla u_t \circ (J-J_0)\rangle \bigg| _{t=0} & = &\\
\label{term31} & \hspace{-8 cm} = & \hspace{-4 cm} \int _{B^{2n}}(-\dive {\xi}) \frac{\partial u^\alpha}{\partial x_{k'}}\frac{\partial u^\alpha}{\partial x_{l'}} [a_{k'l}]([a_{l'k}]+2[\deltat _{l'k}])\\
\label{term32} & \hspace{-8 cm} + & \hspace{-4 cm} 2\int _{B^{2n}} \frac{\partial u^\alpha}{\partial x_{k'}}\frac{\partial u^\alpha}{\partial x_{l''}}\frac{\partial \xi^{l''}}{\partial x_{l'}} [a_{k'l}]([a_{l'k}]+2[\deltat _{l'k}]) \ .
\end{eqnarray}
Therefore we obtain the following equation for $u$:
\begin{eqnarray}
\label{psh} \int _{B^{2n}} \sum _{i,j}(|\nabla u|^2\delta _{ij}-2\langle \frac{\partial u}{\partial x_i},\frac{\partial u}{\partial x_j}\rangle )\frac{\partial \xi ^j}{\partial x_i} & = & \\
\nonumber & \hspace{-6 cm} = & \hspace{-3 cm} -\frac{1}{2} \int _{B^{2n}} \frac{\partial u^\alpha}{\partial x_k}\frac{\partial u^\beta}{\partial x_l}[a_{kl}][b_{\beta\alpha}]\delta _{ij} \frac{\partial \xi ^j}{\partial x_i} \\
\nonumber & \hspace{-6 cm} + & \hspace{-3 cm} \int _{B^{2n}} \frac{\partial u^\alpha}{\partial x_k}\frac{\partial u^\beta}{\partial x_j}[a_{ki}][b_{\beta\alpha}]\frac{\partial \xi^j}{\partial x_i}\\
\nonumber & \hspace{-6 cm} - & \hspace{-3 cm} \frac{1}{2} \int _{B^{2n}} \frac{\partial u^\alpha}{\partial x_{k'}}\frac{\partial u^\alpha}{\partial x_{l'}} [a_{k'l}]([a_{l'k}]+2[\deltat _{l'k}])\delta _{ij}\frac{\partial \xi^j}{\partial x_i} \\
\nonumber & \hspace{-6 cm} + & \hspace{-3 cm} \int _{B^{2n}} \frac{\partial u^\alpha}{\partial x_{k'}}\frac{\partial u^\alpha}{\partial x_j} [a_{k'l}]([a_{ik}]+2[\deltat _{ik}])\frac{\partial \xi^j}{\partial x_i} \ ,
\end{eqnarray}
which is valid for all $\xi\in C^\infty _c(B^{2n},\mathbf{R}^{2n})$. Note that this equation is very similar to the equation for stationary harmonic maps where the right-hand side vanishes, whereas in our case the right-hand side is $O(r)\| \nabla u\| _{L^2}^2$. It is clear that for small enough radii (depending only on $J$) the second order operator involved is only a small perturbation of the Laplacian, and hence elliptic.

\subsubsection{Monotonicity formula for the energy of $u$}

In this subsection we show that from equation (\ref{psh}) we can deduce an almost monotonicity formula for the energy of $u$. We will show that for locally approximable $J$-holomorphic maps $\frac{e^{C\rho}-\rho}{\rho ^{2n-2}}\int _{B_\rho ^{2n}} |\nabla u |^2$ is an increasing function of $\rho$ provided that $\rho$ is smaller than some fixed $r_0$ independent of $u$. Precisely, we will prove the following proposition:\\
\begin{prop}
Let $u: B^{2n}\rightarrow N\subset\mathbf{R}^{2N}$ be a locally approximable $J$-holomorphic map. Then there exist a constant $C>0$ and $r_0\in (0,1]$ both independent of $u$ such that 
$$
\frac{e^{C\tau}-\rho}{\tau ^{2n-2}}\int _{B_\tau ^{2n}} |\nabla u |^2-\frac{e^{C\sigma}-\sigma}{\sigma ^{2n-2}}\int _{B_\sigma ^{2n}} |\nabla u |^2\geq 2\int _{B_\tau\setminus B_\sigma}R^{2-2n}\bigg|\frac{\partial u}{\partial R}\bigg|^2
$$
for all $\sigma$, $\tau \in (0,r_0)$ with $\sigma <\tau$ (here $R=|x|$ and $\frac{\partial u}{\partial R}$ denotes the derivative in the direction of $R$).
\end{prop}

\noindent \textbf{Proof:}\\
The proof will follow the proof for stationary harmonic maps given in \cite{si96}. We use the fact that if $\int a_j\frac{\partial \xi}{\partial x_j}=0$ for all $\xi\in C^\infty _c(B)$, then for a.e. $\rho\in (0,1)$ we have $\int _{B_\rho} a_j \frac{\partial \xi}{\partial x_j}=\int _{\partial B_\rho} \eta \cdot a \xi$. Using this we test the above equation for $u$ (\ref{psh}) with $\xi ^j(x)=x^j$ (so that $\frac{\partial \xi ^j}{\partial x_i}=\delta _{ij})$. This yields the following identity for $u$:
\begin{eqnarray}
\nonumber (2-2n)\rho^{1-2n}\int _{B_\rho} |\nabla u|^2 +\rho ^{2-2n}\int _{\partial B_\rho} |\nabla u |^2 & = & \\
\nonumber & \hspace{-14 cm} = & \hspace{-7 cm} \rho ^{1-2n}\int _{B_\rho} A(x)|\nabla u|^2-\rho ^{2-2n}\int _{\partial B_\rho} B(x)|\nabla u |^2 +2\int _{\partial B_\rho} R^{2-2n}\bigg| \frac{\partial u}{\partial R}\bigg| ^2 \ ,
\end{eqnarray}
where $\int _{B_\rho} A(x)|\nabla u|^2$ and $\int _{\partial B_\rho} B(x)|\nabla u|^2$ denote terms which can be bounded by $C\rho \int _{B_\rho} |\nabla u| ^2$ and $C\rho \int _{\partial B_\rho} |\nabla u|^2$ respectively. This leads to the chain of inequalities below
\begin{eqnarray}
\label{monest} \Bigg| \frac{d}{d\rho}\Bigg(\frac{1}{\rho^{2n-2}}\int _{B_\rho} |\nabla u|^2 - 2\int _{B_\rho} R^{2-2n}\bigg| \frac{\partial u}{\partial R}\bigg| ^2\Bigg) \Bigg| & \leq & \\
\nonumber & \hspace{-8 cm} \leq & \hspace{-4 cm} C\frac{1}{\rho^{2n-2}}\int _{B_\rho} |\nabla u|^2 +C \frac{1}{\rho ^{2n-3}}\int _{\partial B_\rho} |\nabla u|^2 \\
\nonumber & \hspace{-8 cm} \leq & \hspace{-4 cm} C\frac{1}{\rho ^{2n-2}}\int _{B_\rho} |\nabla u|^2 + C\frac{d}{d\rho}\bigg( \rho ^{2n-3} \int  _{B_\rho} |\nabla u|^2 \bigg) \ .
\end{eqnarray}
Using this estimate we conclude that 
\begin{eqnarray}
\nonumber \frac{d}{d\rho}\bigg( \frac{e^{C\rho}}{\rho^{2n-2}}\int _{B_\rho} |\nabla u|^2\bigg) & = & e^{C\rho} \frac{d}{d\rho}\bigg( \frac{1}{\rho^{2n-2}}\int _{B_\rho} |\nabla u|^2\bigg) + \frac{Ce^{C\rho}}{\rho ^{2n-2}}\int _{B_\rho} |\nabla u|^2 \\
\nonumber & \hspace{-4 cm} \geq & \hspace{-2 cm} 2 \frac{d}{d\rho} \Bigg( \int _{B_\rho} R^{2-2n} \bigg|\frac{\partial u}{\partial R}\bigg|^2\Bigg) + \frac{d}{d\rho} \bigg( \frac{1}{\rho^{2n-3}} \int _{B_\rho} |\nabla u|^2\bigg) \ ,
\end{eqnarray}
and integration from $\sigma$ to $\tau$ yields the desired result.\\

\noindent \textbf{Remark:}\\
From equation (\ref{monest}) in the proof of the monotonicity formula we also get the estimate
\begin{eqnarray}
\nonumber \bigg| \frac{d}{d\rho} \bigg( \frac{1}{\rho ^{2n-2}} \int _{B_\rho} |\nabla u|^2 \bigg)\bigg| & \leq & C\frac{1}{\rho ^{2n-2}} \int _{B_\rho} |\nabla u|^2 + C \frac{1}{\rho ^{2n-3}} \int _{\partial B_\rho} |\nabla u|^2 \\
\nonumber & + & 2\int _{\partial B_\rho} R^{2-2n}\bigg| \frac{\partial u}{\partial R}\bigg|^2 \\
\nonumber & \leq & C\frac{1}{\rho ^{2n-2}} \int _{B_\rho} |\nabla u|^2+ C\rho \bigg| \frac{d}{d\rho} \bigg( \frac{1}{\rho ^{2n-2}} \int _{B_\rho} |\nabla u |^2\bigg) \bigg| \\
\nonumber & + & 2\int _{\partial B_\rho} R^{2-2n}\bigg| \frac{\partial u}{\partial R}\bigg|^2 \ .
\end{eqnarray}
Therefore for $\rho$ so small that $\frac{1}{2}\leq 1-C\rho$ we obtain that 
$$
\frac{d}{d\rho} \bigg( \frac{1}{\rho ^{2n-2}} \int _{B_\rho} |\nabla u|^2 \bigg) \leq C\frac{1}{\rho ^{2n-2}} \int _{B_\rho} |\nabla u|^2+ C\int _{\partial B_\rho} R^{2-2n}\bigg| \frac{\partial u}{\partial R}\bigg|^2 \ ,
$$
i.e. that 
$$
\frac{d}{d\rho}\bigg( \frac{e^{-C\rho}}{\rho^{2n-2}} \int _{B_\rho} |\nabla u|^2\bigg) \leq C\int _{\partial B_\rho} R^{2-2n}\bigg| \frac{\partial u}{\partial R}\bigg|^2 \ .
$$
Integrating this inequality from $\sigma$ to $\tau$ we obtain the estimate
\be
\label{altmon} \frac{e^{-C\tau}}{\tau^{2n-2}} \int _{B_\tau} |\nabla u|^2 - \frac{e^{-C\sigma}}{\sigma^{2n-2}} \int _{B_\sigma} |\nabla u|^2 \leq C\int _{B_\tau\setminus B_\sigma} R^{2-2n}\bigg| \frac{\partial u}{\partial R}\bigg|^2 \ ,
\ee
which we will use in the proof below.

\subsubsection{Obtaining a rate of convergence for $u$}
We now give the proof of theorem \ref{J-maps}  in the perturbed situation. As above we use the monotonicity formula to deduce that
$$
C\int _{B_r}R^{2-2n}\bigg|\frac{\partial u}{\partial R}\bigg|^2 \leq \frac{e^{Cr}}{r ^{2n-2}}\int _{B_r} |\nabla u| ^2 -\Theta _u(x_0)< \infty\ .
$$
Contrary to the easy case, this time we cannot slice by $J_0$-holomorphic planes but we have to slice by $J$-holomorphic curves passing through $x_0=0$. The fact that such curves exist for any almost complex structure $J$ was proven in a paper by the second author and G.Tian \cite{riti01} (see their Appendix A). There they also proved that such curves form a singular foliation of  $B^{2n}\subset\mathbf{R}^{2n}$. It is important to note that the tangent plane to such a $J$-holomorphic curve is spanned by two vectors $X$ and $JX$, where near the origin $X$ is only a small perturbation of $\frac{\partial}{\partial R}$. More precisely, we have that any tangent plane is spanned by a vector $X$ for which we have
$$
X=\frac{\partial }{\partial R}+O(r)\frac{\partial}{\partial x_i} \ ,
$$
i.e. for a map $u$ we get that 
$$
u\cdot X=\frac{\partial u}{\partial R}+O(r)\frac{\partial u}{\partial x_i} \ .
$$
For a $J$-holomorphic map we will use the monotonicity formula to show that
$$
\int _{B_r} R^{2-2n} |u \cdot X|^2 <\infty \ .
$$
First note that the remark on $X$ implies that we get the estimate
\begin{eqnarray}
\nonumber \int _0^r d\rho \int _{\partial B_\rho} R^{2-2n}|u\cdot X|^2 & \leq & \int _0 ^r d\rho \int _{\partial B_\rho} \bigg[ R^{2-2n}\bigg| \frac{\partial u}{\partial R}\bigg| ^2 +CR^{3-2n} |\nabla u|^2\bigg] \\
\nonumber & = & \int _{B_\rho} R^{2-2n} \bigg|\frac{\partial u}{\partial R}\bigg|^2 +C \int _{B_\rho} R^{3-2n} |\nabla u|^2 \ .
\end{eqnarray}
The bound on the first term follows directly from the monotonicity formula, whence it remains only to bound the second one. This can be done through integration by parts and applying the monotonicity formula to each term:
\begin{eqnarray}
\nonumber \int _0^r d\rho\ \rho^{3-2n}\int _{\partial B_\rho} |\nabla u|^2 & = & \\
\nonumber & \hspace{-4 cm} = & \hspace{-2 cm} \bigg[\rho^{3-2n}\int _{B_\rho} |\nabla u|^2 \bigg] _0^r+(2n-3)\int _0^r d\rho \ \rho^{2-2n}\int _{B_\rho} |\nabla u|^2 \\
\nonumber & \hspace{-4 cm} \leq & \hspace{-2 cm} Cr \bigg[ \sup_{0<\rho<1}\frac{1}{\rho^{2n-2}}\int _{B_\rho} |\nabla u|^2 \bigg] \\
\nonumber & \hspace{-4 cm} - & \hspace{-2 cm}  \lim _{\rho\rightarrow 0}\bigg(\rho\bigg[ \sup_{0<\rho<1}\frac{1}{\rho^{2n-2}}\int _{B_\rho} |\nabla u|^2 \bigg] \bigg)   <\infty \ .
\end{eqnarray}
In the paper by the second author and G. Tian mentioned above it was shown that the $J$-holomorphic curves foliating $B^{2n}$ are smoothly parametrised by $\mathbf{CP}^{n-1}$. We take $f:B^{2n}\rightarrow \mathbf{CP}^{n-1}$ to be this parametrisation. This map $f$ has the property that its gradient is $O(\frac{1}{r})$; precisely, there are constants $C_1$, $C_2>0$ such that first
$$
C_1-C_2R\leq R^{2n-2} |J_{2n-2}f|\leq C_1+C_2R
$$
(again $J_{2n-2}f$ denotes the $(2n-2)$-Jacobian of $f$) and second $1\leq \frac{C_1+C_2R}{C_1-C_2R}\leq 1+O(r)$. Therefore we can apply the co-area formula to the slicing by $f$ and obtain that 
$$
\int _{\mathbf{CP}^{n-1}} dp \int _{f^{-1}(p)} \frac{R^{2-2n} |u\cdot X|^2}{|J_{2n-2}f|} = \int _{B^{2n}} R^{2-2n} |u\cdot X|^2 <\infty \ ,
$$
which is what we wanted to show. By Fubini's theorem we know that for a.e. $p\in\mathbf{CP}^{n-1}$ we have
$$
\int _{f^{-1}(p)} |u\cdot X|^2 \leq C_2 \int _{f^{-1}(p)} \frac{R^{2-2n} |u\cdot X|^2}{|J_{2n-2} f|} < \infty \ .
$$
Since the metric on the target $g_N$ is compatible with the target complex structure $J_N$ we know that $|u\cdot JX|^2 =|J_N\circ \nabla u(X)|^2=|\nabla u(X)|^2$ and hence from the above argument we also deduce that $\int _{f^{-1}(p)} |u\cdot JX|^2 <\infty$.\\
Next we want to show how the above implies that $\int _{f^{-1}(p)} |\nabla u|^2 <\infty$ for a.e. $p\in\mathbf{CP}^{n-1}$. To see this first note that as $X$ and $JX$ span the tangent space to $f^{-1}(p)$
$$
\nabla u|_{f^{-1}(p)}=(\nabla u\cdot X)X+(\nabla u \cdot JX)JX
$$
from which we conclude that
\be
\label{gradest} |u\cdot X|^2 \leq \frac{1}{2-Cr} |\nabla u|_{f^{-1}(p)}|^2 \leq \frac{2+Cr}{2-Cr} |u\cdot X|^2
\ee
for some $C>0$. Therefore we obtain that for $r$ small enough $\int _{f^{-1}(p)} |\nabla u|^2$ is bounded by some constant times the integrals in the $X$ and $JX$ directions. Hence we have shown that as in the easy case $u|_{f^{-1}(p)}$ is a $J$-holomorphic $W^{1,2}$-map from $(f^{-1}(p), J)$ into $(N, J_N)$. \\
Since $(f^{-1}(p), J)$ is a Riemann surface for each $p\in\mathbf{CP}^{n-1}$ we get that $u|_{f^{-1}(p)}$ is smooth and, since $u|_{f^{-1}(p)}$ satisfies an equation similar to the $\bar{\partial}$-equation, that there exists $\delta >0$ such that for small enough $r$
\be
\int _{f^{-1}(p)\cap B_r} |\nabla u|^2 < \delta \int _{f^{-1}(p)\cap B_{2r}} |\nabla u|^2 \ .
\ee
Together with equation (\ref{gradest}) this gives the following estimate
\bea
\nonumber \int _{f^{-1}(p)\cap B_r} |u\cdot X|^2 & \leq & \frac{1}{2-Cr} \int _{f^{-1}(p)\cap B_r} |\nabla u|^2 < \frac{\delta}{2-Cr} \int _{f^{-1}(p)\cap B_r} |\nabla u|^2 \\
\nonumber & \leq & \delta \frac{2+Cr}{2-Cr} \int _{f^{-1}(p)\cap B_{2r}} |u\cdot X|^2 \ ,
\eea
where we can take $r$ so small that $\delta \frac{2+Cr}{2-Cr}<1$. Integrating this over $\mathbf{CP}^{n-1}$ we get
\be
\int _{\mathbf{CP}^{n-1}} dp \int _{f^{-1}(p)\cap B_r} |u\cdot X|^2 < \deltat \int _{\mathbf{CP}^{n-1}} dp \int _{f^{-1}(p)\cap B_{2r}} |u\cdot X|^2 
\ee
for some $\deltat\in (0,1)$. Applying the co-area formula again we obtain the estimate 
\bea
\int _{B_r} |u\cdot X|^2 |J_{2n-2}f| & = & \int _{\mathbf{CP}^{n-1}} dp \int _{f^{-1}(p)\cap B_r} |u\cdot X|^2\\ 
\nonumber & \hspace{-4 cm}  < & \hspace{-2 cm} \deltat \int _{\mathbf{CP}^{n-1}} dp \int _{f^{-1}(p)\cap B_{2r}} |u\cdot X|^2 = \deltat \int _{B_{2r}} |u\cdot X|^2 |J_{2n-2}f| \ .
\eea
As in the easy case this gives 
$$
\int _{B_r} R^{2-2n} |u\cdot X|^2 < \deltat \int _{B_{2r}} R^{2-2n} |u\cdot X|^2
$$
which implies that there exist $C>0$ and $\gamma \in (0,1]$ such that 
$$
\int _{B_r} R^{2-2n} |u\cdot X|^2\leq C r^\gamma \ .
$$
From this we will deduce a rate of convergence for $\int _{B_r} R^{2-2n} \big|\frac{\partial u}{\partial R}\big|^2$. Note that from $X=\frac{\partial}{\partial R}+O(r)\frac{\partial}{\partial x_i}$ we deduce that 
\bea
\nonumber \int _{B_r} R^{2-2n} \bigg|\frac{\partial u}{\partial R}\bigg|^2 & \leq & \int _{B_r} R^{2-2n} |u\cdot X|^2 \leq Cr^\gamma + C\int _{B_r} R^{3-2n} |\nabla u|^2 \\
\nonumber & \leq & Cr^\gamma +Cr \sup_{0<\rho <1}\frac{1}{\rho^{2n-2}}\int _{B_\rho} |\nabla u|^2 \\
\nonumber & \leq & C(1+\| \nabla u\|^2 _{L^2(B_1)})r^\gamma \ ,
\eea
where for the last estimate we used the monotonicity formlua.\\
To finish the proof we will use the remark at the end of the proof of the monotonicity formula. From estimate (\ref{altmon}) we get that
$$
\frac{e^{-C\tau}}{\tau ^{2n-2}}\int _{B_\tau} |\nabla u|^2 - \Theta _u (x_0) \leq C \int _{B_r} R^{2-2n} \bigg|\frac{\partial u}{\partial R}\bigg|^2 \leq C\Big(1+\| \nabla u\|^2 _{L^2(B_1)}\Big)r^\gamma
$$
which for $\tau$ so that $e^{2C\tau}\leq 2$ becomes
$$
\frac{e^{C\tau}}{\tau ^{2n-2}}\int _{B_\tau} |\nabla u|^2 - e^{2C\tau}\Theta _u (x_0) \leq 2C\Big(1+\| \nabla u\|^2 _{L^2(B_1)}\Big)r^\gamma \ .
$$
Therefore we get that
\bea
\nonumber \frac{e^{C\tau}-\tau}{\tau ^{2n-2}}\int _{B_\tau} |\nabla u|^2 - \Theta _u (x_0) & \leq & 2C\Big(1+\| \nabla u\|^2 _{L^2(B_1)}\Big)r^\gamma + \Big( e^{2C\tau}-1\Big) \Theta _u (x_0) \\
\nonumber & \leq & C\Big(1+\| \nabla u\|^2 _{L^2(B_1)}\Big)r^\gamma + C\| \nabla u\|^2 _{L^2(B_1)} r \\
\nonumber & \leq & C\Big(1+\| \nabla u\|^2 _{L^2(B_1)}\Big)r^\gamma \ ,
\eea
which completes the proof of theorem \ref{J-maps}.

\subsubsection{Uniqueness of tangent maps}

In this subsection we show how theorem \ref{J-maps} implies the uniqueness of tangent maps at all $x_0\in M$.

\begin{prop}\label{mapsunique}
Let $u$, $M$, $N$ be as in theorem \ref{J-maps}. Then $u$ has a unique tangent map at all $x_0\in M$. In particular,
$$
\lim _{r\to 0} u(x_0+rx)
$$
exists weakly in $W^{1,2}$.
\end{prop}

\noindent \textbf{Proof:}\\
The existence of tangent maps for $u$ follows with only minor modifications from a paper by J. Li and G. Tian \cite{liti98}. Now fix $x_0\in M$ and without loss of generality we can assume that we have chosen coordinates so that $x_0=0$. From the proof of theorem \ref{J-maps} we know that there are constants $C_1>0$ and $\gamma >0$ such that
$$
\int _{B_\rho (x_0)} R^{2-2n}\bigg|\frac{\partial u}{\partial R}\bigg| ^2\leq C_1 \rho^\gamma\ .
$$
Setting $u_{0,\rho}(x):=u(\rho x)$ we note that $\int _{B_1}R^{2-2n} \big|\frac{\partial u_{0,\rho}}{\partial R}\big|^2=\int _{B_\rho} R^{2-2n}\big|\frac{\partial u}{\partial R}\big|^2$. Thus
$$
\int _{B_1} R^{2-2n}\bigg|\frac{\partial u_{0,\rho}}{\partial R}\bigg| ^2\leq C_1 \rho^\gamma\ .
$$
From this we deduce that for any $0<\sigma<\tau$ sufficiently small we get
\bea
\nonumber \| u_{0,\tau}-u_{0,\sigma}\|_{L^2(S^{n-1})} & \leq & \int _\sigma ^\tau \bigg\|\frac{\partial u_{0,r}}{\partial r}\bigg\|_{L^2(S^{n-1})} dr\\
\nonumber & \leq & \bigg(  \int _\sigma ^\tau r^{1-\frac{\gamma}{2}} \bigg\|\frac{\partial u_{0,r}}{\partial r}\bigg\|_{L^2(S^{n-1})}^2 dr\bigg)^{\frac{1}{2}} \bigg( \int _\sigma ^\tau r^{-(1-\frac{\gamma}{2})} dr\bigg)^{\frac{1}{2}}\\
\nonumber & \leq & C_2 \tau^{\frac{\gamma}{4}}\ ,
\eea
where the fact that $\int _\sigma ^\tau r^{1-\frac{\gamma}{2}} \big\|\frac{\partial u_{0,r}}{\partial r}\big\|_{L^2(S^{n-1})}^2 dr\leq C_2$ follows from an integration by parts. If we now take a subsequence $\sigma _j\to$ as $j\to\infty$ such that $u_{0,\sigma _j}$ converges to some tangent map $u_\infty$ weakly in $W^{1,2}$ but strongly in $L^2(S^{2n-1})$, then by the triangle inequality and the above estimates for any other sequence $\tau _j\to 0$ we have $u_{0,\tau_j}\to u_\infty$ strongly in $L^2(S^{2n-1})$ and hence that $u_\infty$ is the unique tangent map.


\section{Proofs of theorems \ref{uniqueness} and \ref{rate}}\label{currents}
\subsection{Preliminaries}\label{preliminaries}

We will adopt the standard notation of geometric measure theory (see \cite{fe69}). In this section we will recall some facts about calibrated cycles (of arbitrary dimension and co-dimension in $\Rm$), compute a monotonicity formula for them and investigate the structure of their tangent cones.\\

Since throughout the paper we are assuming that the calibration is at least $C^2$ (but not necessarily closed), we know that for any point $x_0\in\Rm$ there exists a neighborhood depending only on the $C^2$-norm of $\omega$ such that in this neighborhood we can write $\omega (x)=\omega _0(x)+\omega _1(x)$ with $\omega _0(x)=\omega (x_0)$ and $\omega _1(x)=O(|x-x_0|)$.\\

In the special case where $\omega$ is closed and the cycle is integer rectifiable from \cite{fe69} theorem 5.4.3 (2) we know that $C$ satisfies a monotonicity formula. In the case where $C$ is a normal cycle calibrated by a constant calibration, Harvey and Lawson proved a monotonicity formula depending only on this first order information (see \cite{hala82} theorem 5.7). Based on their proof we will now show an almost monotonicity formula in the general case.

\begin{prop}\label{calmon}
Let $C$ be a $p$-dimensional normal cycle in $\mathbf{R}^m$. Assume that $C$ is calibrated by a comass $1$ $p$-form $\omega$. Then there exist $C_1>0$, $r_0>0$ depending only on the $C^2$-norm of $\omega$ such that given $x_0\in\spt C$ for any $0<s<r\leq r_0$ we have
\begin{eqnarray}
\nonumber \frac{e^{C_1r}+C_1r}{r^p}M(C\res B_r(x_0)) & - & \frac{e^{C_1s}+C_1s}{s^p}M(C\res B_s(x_0))\\ 
\nonumber & \hspace{-4 cm} \geq & \hspace{-2 cm} \int _{B_r\setminus B_s(x_0)} \frac{1}{|x-x_0|^p}\sum ^{N(x)}_{i=1} \lambda _i(x)\bigg|\xi _i(x)\wedge \frac{\partial}{\partial r}\bigg|^2 d\|C\|(x) \ ,
\end{eqnarray}
where $C$ is represented by $\int \langle \tau (x), \cdot \rangle d\| C\|$ and $\sum ^{N(x)}_{i=1} \lambda _i(x)\xi _i(x)$ is the decomposition of $\tau (x)$ into a convex sum of calibrated simple vectors.
\end{prop}

\noindent \textbf{Proof:}\\
Using the setting established above we write $\omega (x)=\omega _0(x)+\omega _1(x)$ where $\omega _0(x)=\omega (x_0)$ and $\| \omega _1(x)\| _{C^2}=O(|x-x_0|)$. Note that since $\omega _0$ is a constant $p$-form we know that $\omega _0=\frac{1}{p} d\big(\frac{\partial}{\partial r}\ins \omega _0\big)=\frac{1}{p} L_{\frac{\partial}{\partial r}}\omega _0$ --- here $L_{\frac{\partial}{\partial r}}$ denotes the Lie derivative in the direction of $\frac{\partial}{\partial r}$. Next we take a smooth cut-off function $\phi:\mathbf{R}\rightarrow \mathbf{R}$ such that $\phi (t)=1$ for $t\leq\tfrac{1}{2}$, $\phi (t)=0$ for $t\geq 1$ and $\phi'(t)\leq 0$. Setting $\gamma (x):=\phi \big( \frac{r}{\rho}\big)$ (here and subsequently $r=|x|$), $I(\rho ):=\int _{\Rm} \gamma (x)\langle \tau ,\omega\rangle d\|C\|$ and 
$$
\omega ^t:=\frac{\partial}{\partial r}\ins \bigg(\frac{\partial}{\partial r}\wedge\omega\bigg)
$$
we make the following computation:
\begin{eqnarray}
\nonumber pI(\rho) & = & p\int _{\Rm} \gamma (x)\langle \tau ,\omega _0\rangle d\|C\| + p\int _{\Rm} \gamma (x)\langle \tau ,\omega _1\rangle d\|C\|\\
\nonumber & = & \int _{\Rm}\bigg\langle \tau ,\gamma (x) d\bigg(\frac{\partial}{\partial r}\ins \omega _0\bigg)\bigg\rangle d\|C\| + p\int _{\Rm} \gamma (x)\langle \tau ,\omega _0\rangle d\|C\|\\ 
\nonumber & = & \int _{\Rm}\bigg\langle \tau , d\bigg(\gamma (x) \bigg(\frac{\partial}{\partial r}\ins \omega _0\bigg)\bigg)\bigg\rangle d\|C\|\\
\nonumber & - & \int _{\Rm}\bigg\langle \tau , \frac{d}{dr}(\gamma (x)) dr\wedge d \bigg( \frac{\partial}{\partial r} \ins \omega _0\bigg)\bigg\rangle d\|C\| +p\int _{\Rm} \gamma (x)\langle \tau ,\omega _0\rangle d\|C\|\\
\nonumber  & = & \rho \int _{\Rm} \bigg\langle \tau , \frac{d}{d\rho}\phi \Big( \frac{r}{\rho} \Big)\bigg[ \frac{dr}{r}\wedge\frac{\partial}{\partial r}\ins \omega _0\bigg] \bigg\rangle d\|C\| + p\int _{\Rm} \phi \Big( \frac{r}{\rho} \Big) \langle \tau , \omega _1\rangle d\|C\|\\
\nonumber & = & \rho \int _{\Rm} \bigg\langle \tau , \frac{d}{d\rho}\phi \Big(\frac{r}{\rho}\Big)(\omega _0 - \omega _0^t) \bigg\rangle d\|C\| + p\int _{\Rm} \phi \Big(\frac{r}{\rho}\Big) \langle \tau , \omega _1\rangle d\|C\|\\
\nonumber & = & \rho \int _{\Rm} \bigg\langle \tau , \frac{d}{d\rho}\phi \Big(\frac{r}{\rho}\Big)\omega \bigg\rangle d\|C\| - \rho \int _{\Rm} \bigg\langle \tau , \frac{d}{d\rho}\phi \Big(\frac{r}{\rho}\Big)\omega ^t \bigg\rangle d\|C\| \\
\nonumber & - & \rho \int _{\Rm} \bigg\langle \tau , \frac{d}{d\rho}\phi \Big(\frac{r}{\rho}\Big)\omega _1 \bigg\rangle d\|C\| + \rho \int _{\Rm} \bigg\langle \tau , \frac{d}{d\rho}\phi \Big(\frac{r}{\rho}\Big)\omega _1^t \bigg\rangle d\|C\|\\
\nonumber & + & p\int _{\Rm} \phi \Big(\frac{r}{\rho}\Big) \langle \tau , \omega _1\rangle d\|C\|\ .
\end{eqnarray}
Setting $J(\rho):=\int _{\Rm} \big\langle \tau , \phi \big(\frac{r}{\rho}\big)\omega ^t\big\rangle d\|C\|$ the above computation can be summarised as 
\begin{eqnarray}
\nonumber -\frac{pI(\rho)}{\rho ^{p+1}}+\frac{I'(\rho)}{\rho ^p}-\frac{J'(\rho)}{\rho ^p} & = & \frac{1}{\rho^p}\int _{\Rm} \bigg\langle \tau , \frac{d}{d\rho}\phi \Big(\frac{r}{\rho}\Big)(\omega _1-\omega _1^t) \bigg\rangle d\|C\|\\
\nonumber & - & \frac{p}{\rho ^{p+1}}\int _{\Rm} \phi \Big(\frac{r}{\tau}\Big) \langle \tau , \omega _1\rangle d\|C\|\ .
\end{eqnarray}
If $\rho>0$ is chosen small enough we obtain the following estimate
\begin{eqnarray}
\nonumber \bigg| -\frac{pI(\rho)}{\rho^{p+1}}+\frac{I'(\rho)}{\rho ^p}-\frac{J'(\rho)}{\rho ^p}\bigg| & \leq & \\
\nonumber & \hspace{-4 cm} \leq & \hspace{-2 cm} \frac{C_1}{\rho^p}\int _{\Rm} \bigg|\frac{d}{d\rho}\phi \Big(\frac{r}{\rho}\Big)\bigg|\cdot |x| d\|C\| + \frac{C_2}{\rho^{p+1}}\int _{\Rm} \phi \Big(\frac{r}{\tau}\Big) |x| d\|C\|\\
\nonumber & \hspace{-4 cm} \leq & \hspace{-2 cm} \frac{C_1}{\rho^{p-1}}\int _{\Rm} \frac{d}{d\rho}\phi \Big(\frac{r}{\rho}\Big) d\|C\| + \frac{C_2}{\rho ^p}\int _{\Rm} \phi \Big(\frac{r}{\tau}\Big) d\|C\|\\
\nonumber & \hspace{-4 cm} = & \hspace{-2 cm} \frac{C_1}{\rho^{p-1}} I'(\rho )+\frac{C_2}{\rho^p} I(\rho)\\
\nonumber & \hspace{-4 cm} = & \hspace{-2 cm} \frac{C_3}{\rho^{p-1}} I(\rho) + C_3 \frac{d}{d\rho}\Big(\frac{I(\rho)}{\rho^{p-1}}\Big)\ .
\end{eqnarray}
Using this estimate we have that 
$$
\frac{d}{d\rho}\bigg(\frac{I(\rho)}{\rho ^p}\bigg)+C_3\frac{I(\rho)}{\rho ^p}\geq \frac{1}{\rho ^p}\frac{d}{d\rho}J(\rho)- C_3 \frac{d}{d\rho}\bigg(\frac{I(\rho)}{\rho^{p-1}}\bigg)\ ,
$$
which for $\rho>0$ possibly chosen even smaller becomes
$$
\frac{d}{d\rho}\bigg(\frac{e^{C_3\rho}I(\rho)}{\rho ^p}\bigg)\geq \frac{1}{\rho ^p}\frac{d}{d\rho}J(\rho)- C_3 \frac{d}{d\rho}\bigg(\frac{I(\rho)}{\rho^{p-1}}\bigg)\ .
$$
This implies that for small $\rho$ we have
$$
\frac{d}{d\rho}\bigg(\frac{e^{C_3\rho}+C_3}{\rho ^p}I(\rho)\bigg)\geq \frac{1}{\rho ^p}\frac{d}{d\rho}J(\rho)\ .
$$
Letting $\phi$ increase to the characteristic function of $(-\infty , 1)$, we obtain that the above inequality continues to hold in the sense of distributions, i.e.
$$
\frac{d}{d\rho}\bigg(\frac{e^{C_3\rho}+C_3}{\rho ^p}M(C\res B_\rho (x_0))\bigg)\geq \frac{d}{d\rho}\int _{B_\rho (x_0)}\frac{1}{|x|^p}\langle \tau , \omega ^t\rangle d\|C\|\ .
$$
Since $\tau(x)$ is calibrated by $\omega(x)$ for $\|C\|$-a.e. $x\in B_\rho (x_0)$, with the help of lemma 5.11 in \cite{hala82} we can express the integrand on the right-hand side above as a positive quantity:
$$
\langle \tau (x), \omega ^t (x)\rangle =\sum ^{N(x)}_{i=1} \lambda _i(x)\bigg|\xi _i(x)\wedge \frac{\partial}{\partial r}\bigg|^2\ ,
$$
where $\tau (x)=\sum ^{N(x)}_{i=1} \lambda _i(x)\xi _i(x)$ and the $\xi_j(x)$ are simple vectors calibrated by $\omega (x)$. Integration from $0<s<r\leq r_0$ finishes the proof of the proposition.\\ 

We now look at some implications for tangent cones of calibrated currents. First of all note that the monotonicity formula implies that $\tfrac{1}{r^p}M(C\res B_r(x_0)$ is almost increasing and bounded for small enough radii. If we dilate the current $C$ around $x_0$ by setting 
$$
C_{r,x_0}:=(\lambda _*^{r,x_0} T)\res B_1(x_0)\ ,
$$
where $\lambda _*^{r,x_0}$ means the push-forward by $\lambda ^{r,x_0}(x)=\tfrac{x-x_0}{r}$, then from $M(C_{r,x_0})=r^{-p} M(C\res B_r(x_0))$ we conclude that $M(C_{r,x_0})$ is uniformly bounded as $r$ tends to $0$. From the cycle condition we get that $\partial C_{r,x_0}\res B_1(x_0)=0$ and hence that $N(C_{r,x_0})$ (see \cite{fe69}) is uniformly bounded. Therefore, from the compactness theorem we get that for any sequence of radii $\{r_n\}$ tending to $0$, there exists a subsequence $\{r_{n'}\}$ such that as $n'\rightarrow\infty$
$$
C_{r_{n'},x_0}\xrightarrow{N} C_{\infty, x_0}\ ,
$$
for some normal current $C_{\infty, x_0}$. For reasons that will become apparent later, we call such a limiting current $C_{\infty, x_0}$ a \emph{tangent cone} to $C$ at $x_0$. Note that a priori the limiting object might very well depend on the subsequence chosen and the rest of the paper is devoted to showing that this is not the case, i.e. that there is a unique tangent cone.\\

We now show that the tangent cones to a calibrated current are still calibrated $p$-currents. To see this, we first note that the lower semi-continuity of mass under weak convergence implies that 
\be\label{upmass}
\lim _{r\rightarrow 0}\frac{M(C\res B_r(x_0))}{r^p}=\lim _{r\rightarrow 0}M(C_{r, x_0})\geq M(C_{\infty, x_0})\ .
\ee
However, for currents calibrated by $\omega$ the above inequality can be improved to an equality since $M(C\res B_r(x_0))=C\res B_r(x_0)(\omega)$ which gives that
$$
M(C_{r,x_0})=\frac{1}{r^p} C\res B_r(x_0)(\omega)=C_{r,x_0}(r^p(\lambda ^{r,0x_0})^*\omega)\ ,
$$
i.e. that $C_{r,x_0}$ is calibrated by $r^p(\lambda ^{r,0x_0})^*\omega$. Since $M(C_{r,x_0})\leq C_1<\infty$ we conclude that as $r\rightarrow 0$
$$
|C_{r, x_0}(r^p(\lambda ^{r,0x_0})^*\omega-\omega _0)|\leq C_1\| r^p(\lambda ^{r,0x_0})^*\omega-\omega _0\|_\infty \rightarrow 0
$$
(recall that $\omega$ is in $C^2$). Therefore we obtain that 
\be\label{lowmass}
\lim _{r\rightarrow} M(C_{r, x_0})=\lim_{n'\rightarrow \infty}C_{r_{n'}, x_0}(\omega _0 )=C_{\infty , x_0}(\omega _0)\ .
\ee
Since $\omega_0$ has comass equal to $1$, we conclude that $C_{\infty , x_0}(\omega _0)\leq M(C_{\infty , x_0})$. Combining (\ref{upmass}) and (\ref{lowmass}) we deduce that
$$
\lim _{n'\rightarrow\infty}M(C_{r_{n'}, x_0})=M(C_{\infty , x_0})=C_{\infty , x_0}(\omega _0)\ ,
$$
i.e. that $C_{\infty , x_0}$ is calibrated by $\omega _0$.\\

We continue our discussion of tangent cones by looking at the density of a tangent cone at the origin. From the discussion above we get that 
\bea
\nonumber \frac{1}{r^p}\|C_{\infty , x_0}\|(B_r(0)) & = & \lim _{n'\rightarrow \infty}\frac{1}{r^p}\|C_{r_{n'} , x_0}\|(B_r(x_0))\\
\nonumber & = &  \lim _{n'\rightarrow \infty}\frac{1}{(rr_{n'})^p}\|C\|(B_{rr_{n'}}(x_0))\\
\nonumber & = & \alpha (p)\Theta (\|C\|, x_0)\ ,
\eea
where $\alpha (p)$ denotes the volume of the unit ball in $\mathbf{R}^p$. Thus we conclude that 
$$
\frac{1}{\alpha (p)r^p}\|C_{\infty , x_0}\|(B_r(0))=\Theta (\|C_{\infty , x_0}\|, 0)=\Theta (\|C\|, x_0)\ .
$$
To justify the notion \emph{``tangent cone''} note that currents calibrated by a constant form satisfy a simpler monotonicity for $0<s<r$ (see \cite{hala82} theorem 5.7 for a proof of this, or go through the above proof without the perturbation term):
\bea
\nonumber \frac{1}{r^p}\|C_{\infty , x_0}\|(B_r(0)) - \frac{1}{s^p}\|C_{\infty , x_0}\|(B_s(0)) & & \\ 
\nonumber & \hspace{-8 cm} = & \hspace{-4 cm} \int _{B_r(0)\setminus B_s(0)} \frac{1}{|x|^p}\sum _{j=1}^{N(x)} \lambda _j(x) \bigg| \xi _j (x)\wedge \frac{\partial}{\partial r} \bigg| ^2 d\|C_{\infty , x_0}\|\ ,
\eea
where as above $\tau^\infty (x)=\sum _{j=1}^{N(x)} \lambda _j(x)\xi _j (x)$. From the identity for the density above we conclude that the right-hand side in the monotonicity formula is equal to $0$, i.e. that $\sum _{j=1}^{N(x)} \lambda _j (x) \big| \xi _j (x)\wedge \frac{\partial}{\partial r} \big| ^2=0$ at $\|C_{\infty , x_0}\|$-a.e. $x$. Therefore we get that $\tau ^\infty(x)\wedge \frac{\partial}{\partial r}=0$, $\|C_{\infty , x_0}\|$-a.e., which by the homotopy formula (applied to the affine homotopy from $\lambda ^{1,x_0}$ to $\lambda^{r, x_0}$) implies that $C_{\infty, x_0}$ is a cone.\\

Next we will investigate the support and structure of a tangent cone calibrated by $\omega _0$. To do this we recall the structure theorem for constant $2$-forms of unit comass on $\mathbf{R}^m$ (see \cite{hala82} Theorem 7.16, page 79). For such a $2$-form $\omega _0$ we know that there are coordinates and an almost complex structure $J$, which is compatible with the Euclidean metric, such that $\omega_0$ is the standard symplectic form for this almost complex structure.\\
This implies that for any $x_0\in\spt \|C\|$ we can assume that the coordinates are chosen so that $x_0=0$ and $\omega _0$ is the standard symplectic form on $\mathbf{R}^{2n}\subset\Rm$. Therefore, calibrated $2$-vectors are $0$ in the $\mathbf{R}^{m-2n}$-direction. Using this together with the fact that $\tau^\infty (x)\wedge \frac{\partial}{\partial r}(x)=0$ for $\|C_{\infty, x_0}\|$-a.e. $x\in\Rm$, we deduce that the set of $x\in\spt \|C_{\infty, x_0}\|\cap\mathbf{R}^{n-2m}$ has $\|C_{\infty, x_0}\|$-measure $0$ (since for these $x$, clearly, $\tau^\infty (x)\wedge \frac{\partial}{\partial r}(x)\neq 0$). Thus the support of $\|C_{\infty, x_0}\|$ can be assumed to be contained in $\mathbf{R}^{2n}\subset \Rm$.\\ Furthermore, we know that the approximate tangent planes are $J_0$-holomorphic and thus from $\tau^\infty (x)\wedge \frac{\partial}{\partial r}(x)=0$ we deduce that $\tau^\infty (x)\wedge J_0\frac{\partial}{\partial r}(x)=0$ and hence that $\tau^\infty (x)=\frac{\partial}{\partial r}\wedge J_0\frac{\partial}{\partial r}(x)$ for $\| C_\infty\|$-a.e. $x\in\mathbf{R}^m$. Thus we immediately conclude the following proposition:

\begin{prop}
Let $C_\infty$ be a tangent cone to $C$ at $x_0$ which is calibrated by the $2$-form $\omega _{x_0}$. Then there are coordinates centered at $x_0$ such that for any $2$-form $\psi$, $C_\infty (\psi)$ is of the form:
$$
C_\infty (\psi)=\bigg\langle \phi ; \psi \bigg(\frac{\partial}{\partial r}\wedge J_0\frac{\partial}{\partial r}\bigg)\bigg\rangle\ ,
$$
where $\phi$ is a distribution in $\mathcal{D}'(\mathbf{R}^m)$ with support in $\mathbf{R}^{2n}\subset\mathbf{R}^m$.
\end{prop}

From the fact that $C_\infty$ is a $J_0$-holomorphic cone we can deduce more information on the structure of the above distribution $\phi$:

\begin{prop}\label{conestructure}
Let $\phi$ be the distribution given by the previous proposition. Then $\phi$ is of the form
$$
\langle \phi ;f\rangle = \int _0^1\frac{1}{t}\bigg\langle \Gamma ; \int _{H^{-1}(H(\theta))} f(t, \theta)\ d\Haus ^1(\theta)\bigg\rangle d\Haus ^1(t)\ ,
$$
where $H:S^{2n-1}\to\mathbf{CP}^{n-1}$ is the Hopf fibration and $\Gamma$ is a distribution on $\mathbf{CP}^{n-1}$ determining $\phi$.
\end{prop}

\noindent \textbf{Remark:}\\
Note that the two above propositions combined imply that a tangent cone $C_\infty$ to $C$ at $x_0$ is completely determined by $H_*\partial [C_\infty \res B_1]$, i.e. if two tangent cones $C_\infty ^1$ and $C_\infty ^2$ satisfy $H_*\partial [C^1_\infty \res B_1]=H_*\partial [C^2_\infty \res B_1]$, then $C_\infty ^1=C_\infty ^2$. We will make use of this fact when proving that the tangent cone to $C$ at $x_0$ is unique.\\
In case we have a tangent cone to an integral area-minimizing cycle, F. Morgan proved in \cite{mo82} that then $C_\infty \res B_1$ is a union of $2$-dimensional disks. For $J$-holomorphic integral $1-1$-cycles a more direct proof of this fact was given by the second author and G. Tian in section 2 of \cite{riti03}.\\

\noindent \textbf{Proof of proposition \ref{conestructure}:}\\
From the fact that $C_\infty$ is a cone we know that $\lambda _{{0,r}_*} C_\infty =C_\infty$ for any positive $r$. Thus $\phi$ also satisfies $\langle \phi ; f(t,\theta)\rangle = \langle \phi ; f(rt,\theta)\rangle$ and hence as a distribution $\phi$ is independent of $r$, i.e. $\frac{\partial \phi}{\partial r}=0$. From this one immediately deduces that
$$
\langle \phi ; f(t,\theta)\rangle = \int _0^1 \frac{1}{t}\langle \Sigma ; f(t, \theta)\rangle d\Haus ^1(t)\ ,
$$
where $\Sigma$ is a distribution on $S^{2n-1}\subset \mathbf{R}^{2n}$ and for fixed $t$ we view $f(t,\theta)$ as a function on $S^{2n-1}$.\\
It remains to use the fact that $C_\infty$ is also $J_0$-holomorphic. From this we get that $\langle \Sigma ; f(t,\theta)\rangle =\langle \Sigma ; f(t,J_0\theta)\rangle$ and hence that for any $s\in [0,2\pi]$,
$\langle \Sigma ; f(t,\theta)\rangle =\langle \Sigma ; f(t,e^{is}\cdot\theta)\rangle$, where by $e^{is}\cdot\theta$ we mean the multiplication of each component by $e^{is}$. Thus $\Sigma$ is invariant along the fibers of the Hopf fibration given as $H^{-1}(p)$ for $p\in\mathbf{CP}^{n-1}$. It is then easy to check that $\Sigma$ defines a distribution $\Gamma$ on $\mathbf{CP}^{n-1}$ defined by
$$
\bigg\langle \Gamma ; \int _{H^{-1}(p)} f(r, \tilde{\theta})d\Haus ^1(\tilde{\theta})\bigg\rangle := \bigg\langle \Sigma ; \int _{H^{-1}(H(\theta ))} f(r, \tilde{\theta})d\Haus ^1(\tilde{\theta}) \bigg\rangle
$$
and the proposition holds.


\subsection{Proof of theorem \ref{uniqueness} for calibrated $J_0$-holomorphic $2$-cycles}\label{uniquenessproof}

In this part of the paper we prove theorems \ref{uniqueness} and \ref{rate} in an important special case which we will refer back to when proving the theorems in full generality. The setting we consider now is that $C$ is a $2$-dimensional normal cycle in $\Rmm$ which is calibrated by the standard symplectic form $\omega_0$ of $\Rmm$ --- here we order the coordinates $x_i$ in $\Rmm$ so that for the standard complex structure $J_0$ on $\Rmm$ we get
$$
J_0\cdot \frac{\partial}{\partial x_{2i-1}}=\frac{\partial}{\partial x_{2i}}\ ,\qquad
\omega _0=\sum _{i=1}^m dx_{2i-1}\wedge dx_{2i}\ .
$$
For a calibrated cycle $C$ with tangent vector $\tau$ we thus have that $\tau(x)=\sum _{j=1}^{N(x)} \lambda _j(x)\xi _j (x)$, where the $\xi_j$ are simple $2$-vectors calibrated by $\omega _0$, i.e. where the $\xi_j(x)$ can be written as $\xi_j(x)=v_j(x)\wedge J_0v_j(x)$ (this is an immediate consequence of Wirtinger's inequality, see \cite{fe69} or \cite{hala82}).\\
The proof of the theorems heavily depends on the map $\pi :(\Rmm, J_0)\rightarrow (\mathbf{CP}^{m-1},j_0)$ which we already used for the proof of theorem \ref{J-maps} (recall that $\pi$ is the radial extension of the Hopf fibration). We begin by proving the following lemma:
\begin{lemma}\label{vectest}
Let $\tau$ be in the convex hull of simple vectors calibrated by $\omega_0$ on $\Rmm$, i.e. $\tau=\sum _{j=1}^N \lambda _j\xi _j$ with $\xi _j$ calibrated by $\omega$ and $\sum _{j=1}^N \lambda _j=1$, $0\leq \lambda _j\leq 1$. Then there exists a constant $C_{2m}>0$ depending only on the dimension $2m$ such that for any vector $\zeta\in\Rmm$ we have
$$
\sum _{j=1}^N\lambda |\xi _j\wedge \zeta |^2 \leq \sum _{j=1}^N\lambda |\xi _j\wedge \zeta \wedge J_0\zeta | \leq C_{2m}\sum _{j=1}^N\lambda |\xi _j\wedge \zeta |^2\ .
$$
\end{lemma}

\noindent \textbf{Proof:}\\
Clearly it suffices to prove the lemma when $\tau$ is a single simple vector calibrated by $\omega _0$. In this case we know that $\tau=\xi_1\wedge J_0\xi_1$ for $\xi_1$ of unit length by Wirtinger's inequality (see \cite{fe69}). Since $\xi_1$ and $J_0\xi_1$ are orthonormal we can extend them to an ordered orthonormal basis $\{\xi_1, J_0\xi_1, \xi_2, J_0\xi_2, \ldots , \xi _m, J_0\xi_m\}$ of $\Rmm$. Writing an arbitrary vector $\zeta=\sum _{i=1}^{2m}a_i\xi_i$ in this basis we get that
$$
|\tau\wedge\xi|^2=\Bigg|\sum _{l=1}^{2m}a_l\tau\wedge\xi_l\Bigg|^2=\Bigg|\sum _{l=3}^{2m}a_l\xi_1\wedge J_0\xi_1\wedge\xi_l\Bigg|^2=\sum _{l=3}^{2m}|a_l|^2\ .
$$
Next we compute $|\tau\wedge\zeta\wedge J_0\zeta|$. Writing $\zeta=\sum _{k=1}^{m}(a_{2k-1}\xi_{2k-1}+a_{2k}\xi_{2k})$ and $J_0\zeta = \sum _{l=1}^{m}(a_{2l-1}\xi_{2l}-a_{2l}\xi_{2l-1})$ we obtain
\bea
\nonumber |\tau\wedge\zeta\wedge J_0\zeta| & = & \Bigg| \xi_1\wedge J_0\xi_1\wedge \Bigg[ \sum _{k,l=1}^m (a_{2k-1}a_{2l-1}\xi _{2k-1}\wedge\xi _{2l}-\\
\nonumber & & -(a_{2k-1}a_{2l})\xi _{2k-1}\wedge\xi _{2l-1} + (a_{2k}a_{2l-1})\xi _{2k}\wedge \xi _{2l}-\\
\nonumber & & -(a_{2k}a_{2l})\xi _{2k}\wedge\xi _{2l-1}\Bigg]\Bigg)\\
\nonumber & \leq & C_{2m} \Bigg| \sum _{3\leq l<k}^{2m} a_l^2\xi _1\wedge J_0\xi_1\wedge \xi _l\wedge \xi _k\Bigg|\\
\nonumber & \leq & C_{2m} \sum _{3=l}^{2m} |a_l|^2\ .
\eea
The remaining estimate is done in a similar way.\\

In the next lemma we apply the above result to estimate the mass of $\pi _*[C\res {B_r\setminus B_s(x_0)}]$ for small enough radii.
\begin{lemma}\label{massestimate}
Let $C$ be a normal $2$-cycle in $\Rmm$ which is calibrated by the standard symplectic form $\omega _0$. Then there exist constants $C_1>0$ and $r_0>0$ such that for $0<s<r\leq r_0$ and $x_0\in\spt\|C\|$ we have 
$$
M(\pi _*[C\res {B_r\setminus B_s(x_0)}])\leq C_1\bigg[\frac{M(C\res B_r(x_0))}{r^2}-\frac{M(C\res B_s(x_0))}{s^2}\bigg]\ .
$$
\end{lemma}

\noindent \textbf{Proof:}\\
Without loss of generality we can assume that the coordinates are centered at $x_0$. Let $r_0>0$ be the one obtained in the monotonicity formula. Applying the previous lemma for $\zeta =\frac{\partial}{\partial r}$ we deduce that 
\bea
\nonumber \int _{B_r(0)\setminus B_s(0)}\frac{1}{|x|^2}\sum _{j=1}^{N(x)}\lambda (x) \bigg|\xi _j(x)\wedge \frac{\partial}{\partial r} \bigg|^2 d\|C\| \\
\nonumber & \hspace{-8 cm} \leq & \hspace{-4 cm} \int _{B_r(x)\setminus B_s(x)}\frac{1}{|x|^2}\sum _{j=1}^{N(x)}\lambda (x) \bigg|\xi _j(x)\wedge \frac{\partial}{\partial r} \wedge J_0\frac{\partial}{\partial r} \bigg| d\|C\| \\
\nonumber & \hspace{-8 cm} \leq & \hspace{-4 cm} C_{2m} \int _{B_r(x)\setminus B_s(x)}\frac{1}{|x|^2} \sum _{j=1}^{N(x)}\lambda (x) \bigg|\xi _j(x)\wedge \frac{\partial}{\partial r} \bigg|^2 d\|C\|\ .
\eea
Combining this with the special case of the monotonicity formula we obtain that 
\bea
\nonumber \int _{B_r(x)\setminus B_s(x)}\frac{1}{|x|^2}\sum _{j=1}^{N(x)}\lambda (x) \bigg|\xi _j(x)\wedge \frac{\partial}{\partial r} \wedge J_0\frac{\partial}{\partial r} \bigg| d\|C\| \\
\nonumber & \hspace{-8 cm} \leq & \hspace{-4 cm} C_{2m} \bigg[ \frac{M(C\res B_r(0))}{r^2}-\frac{M(C\res B_s(0))}{s^2} \bigg]\ .
\eea
Since $\pi$ is close to $H\circ\frac{x}{|x|}$ in $C^2$-norm we deduce that 
$$
\langle \tau (x), \bigwedge ^2\pi(x)\rangle =\frac{1}{|x|^2}\bigg|\tau(x)\wedge \frac{\partial}{\partial r} \wedge J_0\frac{\partial}{\partial r} \bigg|
$$
and the proof of the lemma is completed.\\

We now show how we can prove uniqueness of tangent cones (i.e. theorem \ref{uniqueness}) in this special case with the help of these two lemmas. Note that since $r^{-2}M(C\res B_r(x_0))$ is increasing in $r$, the above lemma implies that
\be\label{massprojection}
\lim _{r\to 0}\lim _{s\to 0} M(\pi _*[C\res {B_r\setminus B_s(x_0)}])=0\ .
\ee
Now suppose that $\{r_i\}$ and $\{s_i\}$ are two sequences converging to $0$, where we can assume that $s_i<r_i\leq r_0$ for all $i$. Also assume that the sequences are chosen so that
$$
C_{r_i,x_0}\xrightarrow[i\to\infty]{\mathbf{F}} C^1_{\infty, x_0}\quad \text{and}\quad C_{s_i,x_0}\xrightarrow[i\to\infty]{\mathbf{F}} C^2_{\infty, x_0}\ .
$$
We now want to show that $C^1_{\infty, x_0}=C^2_{\infty, x_0}$. To see this first note that the map $\pi$ can be seen as a composition of the maps $\frac{x}{|x|}:\Rmm\to S^{2m-1}$ and $H:S^{2m-1}\to\mathbf{CP}^{m-1}$ (the Hopf map), i.e. $\pi=H\circ \frac{x}{|x|}$. Then (\ref{massprojection}) implies that
$$
\lim _{i\to\infty} M\bigg(H _*\bigg(\frac{x}{|x|}_*[C\res {B_{r_i}\setminus B_{s_i}(x_0)}]\bigg)\bigg)=0\ ,
$$
which immediately gives
$$
H _*\bigg(\frac{x}{|x|}_*[C\res {B_r\setminus B_s(x_0)}]\bigg)\xrightarrow[i\to\infty]{w}0
$$
and hence that 
\be\label{weakprojection}
H _*\bigg(\partial\frac{x}{|x|}_*[C\res {B_r\setminus B_s(x_0)}]\bigg)\xrightarrow[i\to\infty]{w}0\ .
\ee
Since for almost all $0<s<r$ we have
$$
\partial \frac{x}{|x|}_*[C\res {B_r\setminus B_s(x_0)}] = \partial \lambda _*^{r,x_0} C\res B_r(x_0) - \partial \lambda _*^{s,x_0} C\res B_s(x_0)\ ,
$$
from line (\ref{weakprojection}) we deduce that
$$
H_*\bigg(\partial \lambda _*^{r,x_0} C\res B_r(x_0) - \partial \lambda _*^{s,x_0} C\res B_s(x_0)\bigg)\xrightarrow[i\to\infty]{w}0\ .
$$
From the choice of sequences $\{r_i\}$ and $\{s_i\}$ it now follows that for the tangent cones $C_{\infty, x_0}^1$ and $C_{\infty, x_0}^2$we have
$$
H_*\big( \partial C_{\infty, x_0}^1\big)=H_*\big( \partial C_{\infty, x_0}^2\big)\ .
$$
For integer rectifiable currents the result now follows immediately, since the boundaries of their tangent cones are unions of great circles contained in $J$-holomorphic planes (see the remark after proposition \ref{conestructure}). For normal currents this follows from proposition \ref{conestructure}. Thus the proof of theorem \ref{uniqueness} is completed for this special case.\\


\subsection{Obtaining a rate for integer rectifiable $2$-cycles calibrated by $\omega _0$}\label{proofrate}

The remaining part of this chapter is devoted to the proof of theorem \ref{rate} for integer rectifiable $2$-cycles calibrated by $\omega _0$. From lemma \ref{massestimate} above we know that it suffices to prove a rate for $M(\pi _*C\res B_r(x_0))$ since
\bea
\nonumber \frac{M(C\res B_r(x_0))}{r^2}-\alpha (2)\Theta (\|C\|, x_0) & = & \int _{B_r(x_0)}\frac{1}{|x-x_0|^2} \bigg| \tau\wedge\frac{\partial}{\partial r}\bigg|^2 d\|C\| \\
\nonumber & \leq & M(\pi _*C\res B_r(x_0))\ .
\eea
From now on, by $\omega$ we will denote the standard symplectic $2$-form which is compatible with the metric on $(\mathbf{CP}^{m-1}, j_0)$. We will deduce the rate for $M(\pi _*C\res B_r(x_0))$ from a rate for 
$$
\int _{B_r(x_0)} \pi^*\omega_{|C}  d\|C\| \ .
$$
First note that for all $0<r\leq r_0$
\bea
\nonumber \int _{B_r(x_0)} \pi^*\omega_{|C}  d\|C\| & = & \pi _*[C\res B_r(x_0)](\omega)\\
\nonumber & \leq & \|\omega\| M(\pi _* C\res B_r(x_0))\\
\nonumber & \leq & C_1\bigg[\frac{M(C\res B_{r_0}(x_0))}{{r_0}^2}-\alpha (2)\Theta (\|C\|, x_0)\bigg]<+\infty\ .
\eea
Furthermore, since the map $\pi$ is $J_0$-$j_0$-holomorphic and $C$ is $J_0$ holomorphic, we know that pointwise ($\|C\|$-a.e.) we have $|\nabla \pi_{|C}|^2(x)=\pi^*\omega_{|C}(x)$, i.e. that for all $0<r\leq r_0$ we have 
$$
\int _{B_r(x_0)} |\nabla \pi_{|C}|^2 d\|C\|\leq C_1<+\infty\ .
$$

As mentioned in the remark at the end of section \ref{preliminaries} \cite{riti03} showed that the tangent cones to $C$ at $x_0$ are unions of $Q$ $J_0$-holomorphic disks (here $Q=\Theta(\|C\|, x_0)$). Arguing like in the proof of lemma III.1 (the part to obtain equation III.3) in \cite{riti03} from our theorem \ref{uniqueness} one can deduce the following lemma proved in \cite{riti03}:
\begin{lemma}\label{cone}
Let $C$ be an integer rectifiable $2$-cycle calibrated by $\omega _0$. Let $x_0\in\spt \|C\|$ and let $D_1,\ldots ,D_Q$ be $2$-disks calibrated by $\omega _0$ such that $C_{\infty , x_0}=\bigoplus _{i=1}^QD_i$. Then given $\epsilon >0$ there exists $\rho_\epsilon >0$ such that for any $0<\rho \leq\rho _\epsilon$ and any $\psi\in C_0^\infty \big(\bigwedge ^2 \big( B^{2m}_1\setminus \{x\in B^{2m}_1\ :\ \dist (x,\cup _i D_i)\leq \epsilon |x|\} \big) \big)$ we have 
$$
C_{\rho ,x_0}(\psi)=0\ .
$$
\end{lemma}

\noindent{\bf Remark:}\\
It is important to note that the proof of this lemma in \cite{riti03} only depends on uniqueness of tangent cones, the structure of tangent cones and the monotonicity formula and not on the epiperimetric inequality by B. White \cite{wh83}. We will give a proof of this lemma in a more general setting in the last section of this paper.\\

From the lemma it immediately follows that there exists a cone $K$ centered at $x_0$ such that for $r>0$ small enough we have
$$
M(\pi _* [C\res B_r(x_0)\cap K])=0
$$
and hence there exists a small ball $B$ in $\mathbf{CP}^{m-1}$ such that $B\subset\subset \pi(K)$ and $(\pi _*[T\res B_r(x_0)])\res B=0$. Therefore we can find a smooth $1$-form $\alpha$ on $\mathbf{CP}^{m-1}$ such that $\omega =d\alpha$ on $\mathbf{CP}^{m-1}\setminus B$ --- recall that $\omega$ generates $H^2(\mathbf{CP}^{m-1};\mathbf{R})$, i.e. that the above steps are necessary to obtain the existence of $\alpha$.\\

Before proving the theorem we need one more lemma:
\begin{lemma}\label{goodslice}
For $C$ as in the previous lemma, let $\omega=d\alpha$ on $\mathbf{CP}^{m-1}\setminus B$. Given $0<r\leq r_0$, there exists $\rho _0\in \big[\tfrac{r}{2}, r\big]$ such that the following hold:
\begin{enumerate}
\item $M(\langle C, |\cdot |, \rho_0\rangle )\leq C_1\rho _0$\ ;
\item $\int _{\Rmm} |\nabla \pi _{|C}|^2 d\langle C, |\cdot |, \rho_0\rangle \leq \frac{1}{\rho _0} \int _{B_r(x_0)\setminus B_{\frac{r}{2}}(x_0)} |\nabla \pi _{|C}|^2 d\|C\|$\ ,
\end{enumerate}
where $\langle C, |\cdot |, \rho_0\rangle$ denotes the slice current of $C$ (see \cite{fe69} chapter 4.3).
\end{lemma}

\noindent{\bf Proof:}\\
From \cite{fe69} 4.2.1 and the monotonicity formula we know that there exists $C_1>0$ such that 
\bea
\nonumber \int _\frac{r}{2}^r \frac{1}{\rho}M(\langle C, |\cdot |, \rho \rangle)d\rho & \leq & \frac{2}{r}M(C\res B_r\setminus B_{\frac{r}{2}}) \leq \frac{2}{r} M(C\res B_r)\\
\nonumber & \leq & C_1\frac{r}{2}\ ,
\eea
i.e.
$$
\avint _{\!\!\! \frac{r}{2}}^{\, r} \frac{1}{\rho}M(\langle C, |\cdot |, \rho \rangle)d\rho \leq C_1\ .
$$
Also, setting $\phi (r):=\int _{B_r(x_0)\setminus B_{\frac{r}{2}}(x_0)} |\nabla \pi _{|C}|^2d\|C\|$ we get that 
$$
\int _{\frac{r}{2}}^r\frac{r}{\phi (r)}\int _{\Rmm} |\nabla \pi _{|C}|^2 d\langle C, |\cdot |, \rho \rangle d\rho \leq \frac{r}{\phi (r)} \phi (r)= r\ .
$$
Therefore we obtain an estimate for the following average integral:
$$
\avint _{\!\!\! \frac{r}{2}}^{\, r} \bigg[ \frac{1}{\rho}M(\langle C, |\cdot |, \rho \rangle)+\frac{r}{\phi (r)} \int _{\Rmm} |\nabla \pi _{|C}|^2 d\langle C, |\cdot |, \rho \rangle \bigg]d\rho \leq C_2\ ,
$$
which implies that there exists $\rho _0\in\big[\tfrac{r}{2}, r\big]$ with 
$$
\frac{1}{\rho_0}M(\langle C, |\cdot |, \rho_0 \rangle)+\frac{r}{\phi (r)} \int _{\Rmm} |\nabla \pi _{|C}|^2 d\langle C, |\cdot |, \rho_0 \rangle \leq C_2
$$
and the lemma follows immediately.\\

Now we are in position to prove the theorem. From the fact that $\omega =d\alpha$ on the support of $\pi _*[C\res B_r(x_0)]$, $\partial C=0$, the construction of $\rho _0$ above and the Cauchy-Schwarz inequality we derive the following estimate:
\bea
\nonumber \int _{B_{\frac{r}{2}}(x_0)} |\nabla \pi _{|C}|^2 d\|C\| & \leq & \int _{B_{\rho _0}(x_0)} |\nabla \pi _{|C}|^2 d\|C\| = \int _{B_{\rho _0}(x_0)} \pi^*\omega _{|C} d\|C\| \\
\nonumber & \hspace{-7 cm} = & \hspace{-3.5 cm}  \int _{\Rmm} \pi^*\alpha _{|C} d\langle C, |\cdot |, \rho_0 \rangle \\
\nonumber & \hspace{-7 cm} \leq & \hspace{-3.5 cm} C_1\int _{\Rmm} |\pi_{|C}-\bar{\pi}_{|C}|\ |\nabla \pi_{|C}| d\langle C, |\cdot |, \rho_0 \rangle\\
\nonumber & \hspace{-7 cm} \leq & \hspace{-3.5 cm} C_1\bigg(\int _{\Rmm} |\pi_{|C}-\bar{\pi}_{|C}|^2 d\langle C, |\cdot |, \rho_0 \rangle\bigg)^{\frac{1}{2}} \bigg( \int _{\Rmm} |\nabla \pi_{|C}|^2 d\langle C, |\cdot |, \rho_0 \rangle\bigg)^{\frac{1}{2}}\ .
\eea
Having obtained this inequality we now use the fact that $\langle C, |\cdot |, \rho_0 \rangle$ is an integer rectifiable cycle. From the decomposition theorem for $1$-dimensional integer rectifiable cycles (see \cite{fe69} 4.2.25) we deduce that one can write $\langle C, |\cdot |, \rho_0 \rangle$ as $\langle C, |\cdot |, \rho_0 \rangle =\sum _{i=1}^\infty T_i$ where each $T_i$ is an indecomposable $1$-cycle and $M(\langle C, |\cdot |, \rho_0 \rangle)=\sum _{i=1}^\infty M(T_i)$. Furthermore, for each $i\in\mathbf{N}$ there exists $f_i:\mathbf{R}\to\Rmm$ with $\Lip (f_i)\leq 1$ such that $T_i={f_i}_*[[0,M(T_i)]]$. Using these facts one can apply Poincar\'e's inequality on each of the pieces, where because of $\Lip (f_i)\leq 1$ we can take the constants in the estimate equal to $1$ independent of $i$. Thus we get
\bea
\nonumber \bigg(\int _{\Rmm} |\pi_{|C}-\bar{\pi}_{|C}|^2 d\langle C, |\cdot |, \rho_0 \rangle\bigg)^{\frac{1}{2}} & = & \bigg( \sum _{i=0}^\infty \int _{\Rmm} |\pi_{|C}-\bar{\pi}_{|C}|^2 d\|T_i\| \bigg)^{\frac{1}{2}} \\
\nonumber & \hspace{-4 cm} \leq & \hspace{-2 cm} \bigg( \sum _{i=0}^\infty M(T_i)^2 \int _{\Rmm} |\nabla \pi_{|C}|^2 d\|T_i\| \bigg)^{\frac{1}{2}}\\
\nonumber & \hspace{-4 cm} \leq & \hspace{-2 cm} M(\langle C, |\cdot |, \rho_0 \rangle)\ \bigg( \int _{\Rmm} |\nabla \pi_{|C}|^2 d\langle C, |\cdot |, \rho_0 \rangle\bigg)^{\frac{1}{2}}\ .
\eea
Combining the above inequalities with lemma \ref{goodslice} we get that
\bea
\nonumber \int _{B_{\frac{r}{2}}(x_0)} |\nabla \pi _{|C}|^2 d\|C\| & \leq & M(\langle C, |\cdot |, \rho_0 \rangle)\ \int _{\Rmm} |\nabla \pi_{|C}|^2 d\langle C, |\cdot |, \rho_0 \rangle \\
\nonumber & \leq & C_1 \rho _0\ \frac{1}{\rho _0} \int _{B_r(x_0)\setminus B_{\frac{r}{2}}(x_0)} |\nabla \pi _{|C}|^2 d\|C\|
\eea
which immediately implies that
$$
\int _{B_{\frac{r}{2}}(x_0)} |\nabla \pi _{|C}|^2 d\|C\| \leq \frac{C}{C+1} \int _{B_r(x_0)} |\nabla \pi _{|C}|^2 d\|C\| \ .
$$
Now a standard iteration argument (see the book by M. Giaquinta \cite{giaq83} for details) shows that there exist $C_1>0$ and $\gamma \in (0,1]$ with
$$
\int _{B_r(x_0)} \pi^*\omega _{|C} d\|C\|=\int _{B_r(x_0)} |\nabla \pi _{|C}|^2 d\|C\|\leq C_1r^\gamma \ .
$$
Recall that $\pi$ is $J_0$-$j_0$-holomorphic and hence that $\pi _*[C\res B_r(x_0)]$ is a $j_0$-holomorphic $2$-current (although not a cycle). Therefore $\pi _*[C\res B_r(x_0)]$ is calibrated by $\omega$ so that the above estimate implies
$$
M(\pi _*[C\res B_r(x_0)])=\pi _*[C\res B_r(x_0)](\omega)=\int _{B_r(x_0)} \pi^*\omega _{|C} d\|C\|\leq  C_1r^\gamma \ ,
$$
which together with the monotonicity formula completes the proof of theorem \ref{rate} in this special case.\\

\noindent{\bf Remark:}\\
There are two steps in the above proof which (for the time being) cannot be carried out for calibrated normal cycles in general. The one is lemma \ref{cone}, the other the decomposition theorem for integer rectifiable $1$-cycles which we used to show that Poincar\'e's inequality was valid. It would be interesting to know whether one can nevertheless still use the above strategy to prove a rate of convergence.


\subsection{Uniqueness of tangent cones for calibrated $2$-cycles}

The aim of this section is to prove theorem \ref{uniqueness} in the general case. The approach will be similar to the one given in the previous section.\\
The setting we work in is described as follows: first fix $x_0\in\spt \|C\|$ and assume that we have chosen coordinates so that $\omega _0=\omega (x_0)$ is the standard symplectic form on $\Rnn\subset\Rm$ (see the end of section \ref{preliminaries}). Furthermore, we can assume that we work in a small ball $B_{r_0}(0)$ so that on this ball we have $\| \omega (x)-\omega_0 \| _{C^2}=O(|x|)$ and the almost monotonicity formula holds true.\\

Using this setting we prove the following lemma:
\begin{lemma}\label{almostcal}
If $x\in\spt\|C\|$ and if $\tau (x)$ is calibrated by $\omega (x)=\omega _0+\omega _1(x)$, with $\omega _0$ as above, then we can write $\tau (x)$ as $\tau _0(x)+\tau _1(x)$ so that $\| \tau_0(x)\|=1$, $\omega _0(\tau _0(x))=1$ and $\| \tau _1(x)\|=O(|x|^{\frac{1}{2}})$.
\end{lemma}

\noindent{\bf Proof:}\\
Since $\tau (x)$ is calibrated by $\omega (x)$ we can write $\tau (x)=\sum _{j=1}\lambda _j(x) \xi _j(x)\wedge \frac{\partial}{\partial r}$, where each $\xi _j(x)$ is again calibrated by $\omega (x)$. To show the lemma it therefore clearly suffices to show the lemma for calibrated simple vectors $\xi (x)$. From the fact that $\omega (x)$ calibrates $\xi (x)$ and $\omega (x)=\omega _0 +\omega _1(x)$ we deduce that $\omega _0(\xi (x))=1+O(|x|)$. To construct a simple $2$-vector $\xi _0$ close to $\xi$ we first orthogonally project $\xi (x)$ onto $\Rnn\subset \Rm$ to obtain a simple $2$-vector $\tilde{\xi}(x)$. Then we get that $\omega _0(\tilde{\xi}(x))=1+O(|x|)$ and denoting the projection of $\xi (x)$ onto ${\Rnn}^\perp$ by $\tilde{\xi}^\perp$ we know that therefore $\| \xi (x)-\tilde{\xi}(x)\|=O\Big(|x|^\frac{1}{2}\Big)$. Setting $\bar{\xi}(x):=\frac{\tilde{\xi}}{\| \tilde{\xi}\|}$ we know that there are orthonormal vectors $v(x)$, $w(x)\in \Rnn\subset\Rm$ such that $\bar{\xi}(x)=v(x)\wedge w(x)$. Then the $2$-vector $\xi _0 (x):=v(x)\wedge J_0 v(x)$ is of mass $1$ and calibrated by $\omega _0$. Thus it remains to show that $\|\xi (x)-\xi _0\|=O\Big( |x|^{\frac{1}{2}}\Big)$. To see this note that it suffices to estimate $\|\tilde{\xi} (x)-\xi _0\|^2=\| [\|\tilde{\xi}\| w(x)-J_0 v(x)]\|^2$. This follows from $\omega _0(\tilde {\xi}(x))=1+O(|x|)$ where we deduce that $\| [\|\tilde{\xi}\| w(x)-J_0 v(x)]\|^2=O(|x|)$ and the proof of the lemma is completed.\\

Combining this lemma with lemma \ref{vectest} we immediately obtain the following estimate:
\begin{cor}
For $x$ and $\tau (x)=\tau _0(x)+\tau _1(x)$ as in the lemma above there exists $C_1>0$ such that for any vector $\zeta \in \Rnn\subset\Rm$ we have
\bea
\nonumber |\tau\wedge \zeta |^2 & \leq & 2|\tau_0\wedge \zeta |^2+2|\tau _1\wedge \zeta |^2\\
\nonumber & \leq & 2|\tau _0\wedge \zeta \wedge J_0\zeta | + 2|\tau _1\wedge \zeta |^2\\
\nonumber & \leq & C_1 |\tau _0\wedge \zeta |^2 + 2|\tau _1\wedge \zeta |^2\\
\label{wedgeest} & \leq & C_1 |\tau\wedge\zeta|^2 + (C_1+2) |\tau _1\wedge \zeta |^2\ .
\eea
\end{cor}

As in the easy case, from line \ref{wedgeest} in the corollary we get the following estimate:
\bea
\nonumber \int _{B_r(x)\setminus B_s(x)}\frac{1}{|x|^2}\sum _{j=1}^{N(x)}\lambda (x) \bigg|\xi ^0_j(x)\wedge \frac{\partial}{\partial r} \wedge J_0\frac{\partial}{\partial r} \bigg| d\|C\| & &\\
\nonumber & \hspace{-16 cm} \leq & \hspace{-8 cm} C_{2m} \int _{B_r(x)\setminus B_s(x)}\frac{1}{|x|^2} \sum _{j=1}^{N(x)}\lambda (x) \bigg[\bigg|\xi _j(x)\wedge \frac{\partial}{\partial r} \bigg|^2 + \bigg|\xi _j^1(x)\wedge \frac{\partial}{\partial r} \bigg|^2\bigg] d\|C\|\ .
\eea
The first term on the right-hand side can again be estimated by the almost monotonicity formula. For the term on the right-hand side we will use lemma \ref{almostcal} to show:
\begin{lemma}\label{pertlemma}
For $C$ and $\tau ^1(x)$ as above there exist $C_1>0$ and $r_0>0$ such that $0<s<r\leq r_0$ implies
$$
\int _{B_r(x)\setminus B_s(x)}\frac{1}{|x|^2} \sum _{j=1}^{N(x)}\lambda (x) \bigg|\xi _j^1(x)\wedge \frac{\partial}{\partial r} \bigg|^2d\|C\| \leq C_1 r\ .
$$
\end{lemma}

\noindent{\bf Proof:}\\
First note that by lemma \ref{almostcal} there exists $C_1>0$ with
$$
\int _{B_r(x)\setminus B_s(x)}\frac{1}{|x|^2} \sum _{j=1}^{N(x)}\lambda (x) \bigg|\xi _j^1(x)\wedge \frac{\partial}{\partial r} \bigg|^2d\|C\| \leq C_1 \int _{B_r(x)\setminus B_s(x)}\frac{1}{|x|} d\|C\|\ ,
$$
for $0<s<r\leq r_0$. From 2.5.18 (3) in \cite{fe69} and the almost monotonicity formula we get the following inequalities:
\bea
\nonumber \int _{B_r(x)\setminus B_s(x)}\frac{1}{|x|} d\|C\| & \leq & \int _s^r \frac{1}{\rho} \frac{d}{d\rho}\big[ \|C\|(B_\rho(0))\big]\, d\rho\\
\nonumber & \leq & \frac{1}{\rho} \| C\|(B_\rho(0))\bigg|_s^r + \int _s^r \frac{1}{\rho^2} \big[\|C\| (B_\rho (0))\big]\, d\rho\\
\nonumber & \leq & \rho \sup _{0<\rho\leq r_0}\frac{\| C\|(B_\rho(0))}{\rho ^2} \bigg|_s^r + \sup _{0<\rho\leq r_0}\frac{\| C\|(B_\rho(0))}{\rho ^2} (r-s)\\
\nonumber & \leq & C_1' (r-s)\ ,
\eea
establishing the lemma.\\

Therefore, we again deduce that
$$
\lim _{r\to 0}\lim _{s\to 0} M(\pi _*[C\res {B_r\setminus B_s(x_0)}])=0\ ,
$$
where this time we extended the map $\pi$ to $\pi :\Rnn\times\mathbf{R}^{m-2n}\to\mathbf{CP}{n-1}\times \mathbf{R}^{m-2n}$ by sending $(x_1, \ldots , x_{2n}, x_{2n+1}, \ldots , x_m)$ to $(\pi(x_1,\ldots , x_{2n}), 0)$. From the last part of section \ref{preliminaries} we know that tangent cones to $C$ are entirely determined by their $\Rnn$-components and hence the uniqueness of tangent cones follows by exactly the same arguments as in section \ref{uniquenessproof}.


\subsection{Proof of theorem \ref{rate}}

Using arguments similar to the ones from the previous section we now prove theorem \ref{rate} in the full generality as stated. In section \ref{proofrate} we used the exact monotonicity formula to deduce a rate of convergence from a rate of convergence for $M(\pi _* C\res B_r(x_0))$. In the general case we have to use an estimate like estimate (\ref{altmon}) in the perturbation argument for $J$-holomorphic maps, i.e. there exist $C_1>0$ and $r_0>0$ with
\bea
\nonumber \frac{e^{-C_1r}}{r^2}M(C\res B_r(x_0)) & - & \frac{e^{-C_1s}}{s^2}M(C\res B_s(x_0))\\
\label{comparerate} & \leq & C_1 \int _{B_r(x_0)}\frac{1}{|x-x_0|^2} \bigg| \tau\wedge\frac{\partial}{\partial r}\bigg|^2 d\|C\| \ .
\eea
This estimate follows from the proof of the monotonicity formula in the same fashion as estimate (\ref{altmon}) for $J$-holomorphic maps. Using lemmas \ref{almostcal} and \ref{pertlemma} we can see that
$$
\int _{B_r(x_0)}\frac{1}{|x-x_0|^2} \bigg| \tau\wedge\frac{\partial}{\partial r}\bigg|^2 d\|C\|  \leq M(\pi _* C\res B_r(x_0)) + C_1\, r \ ,
$$
and hence that a rate of convergence would still follow from a rate of convergence for $M(\pi _* C\res B_r(x_0))$.\\
Since lemma \ref{goodslice} is still valid in this case (its proof depended only on the monotonicity formula and general estimates from \cite{fe69}), it remains only to show two steps. The one is to prove lemma \ref{cone} in this context, the other to show that the fact that $C$ is not exactly $J_0$-holomorphic near $x_0$ does not matter. We begin by deducing the latter from lemma \ref{gencone} stated below; in fact, denoting $C$ with orientation $\tau_0(x)$ by $C_0$, we have
\bea
\nonumber \int _{B_{\frac{r}{2}}(x_0)} |\nabla \pi _{|C}|^2 d\|C\| & \leq & \int _{B_{\rho _0}(x_0)} |\nabla \pi _{|C}|^2 d\|C\|\\
\nonumber & \leq & \int _{B_{\rho _0}(x_0)} \pi^*\omega _{|C_0} d\|C\| + C_1'\rho _0 \int _{B_{\rho _0}(x_0)} |\nabla \pi _{|C}|^2 d\|C\|\\
\nonumber & \leq & \int _{B_{\rho _0}(x_0)} d\pi^*\alpha _{|C_0} d\|C\| + C_1'r \int _{B_r (x_0)} |\nabla \pi _{|C}|^2 d\|C\|\\
\nonumber & \leq & \int _{B_{\rho _0}(x_0)} d\pi^*\alpha _{|C} d\|C\| + C_1r \int _{B_r (x_0)} |\nabla \pi _{|C}|^2 d\|C\|\ .
\eea
Now the same strategy as in section \ref{proofrate} applies to show
$$
\int _{B_{\frac{r}{2}(x_0)}} |\nabla \pi _{|C}|^2\, d\|C\| \leq \frac{C_1+C_1 r}{C_1+1} \int _{B_r(x_0)} |\nabla \pi _{|C}|^2\, d\|C\|\ ,
$$
which (having chosen $r$ small enough) implies the theorem the same way as before.\\
To complete the proof of theorem \ref{rate} it remains to show the following lemma which is a part of lemma III.1 in \cite{riti03}. Their proof, although given for calibrated $J$-holomorphic currents, also works in the general case and we include it for the sake of completeness.
\begin{lemma}\label{gencone}
Let $C$ be an integer rectifiable $2$-cycle calibrated by $\omega $. Let $x_0\in\spt \|C\|$ and let $D_1,\ldots ,D_Q$ be $2$-disks calibrated by $\omega _0$ such that $C_{\infty , x_0}=\bigoplus _{i=1}^QD_i$. Then given $\epsilon >0$ there exists $\rho_\epsilon >0$ such that for any $0<\rho \leq\rho _\epsilon$ and any $\psi\in C_0^\infty \big(\bigwedge ^2 \big( B^m_1\setminus \{x\in B^m_1\ :\ \dist (x,\cup _i D_i)\leq \epsilon |x|\} \big) \big)$ we have 
$$
C_{\rho ,x_0}(\psi)=0\ .
$$
\end{lemma}
\noindent{\bf Proof:}\\
For an argument by contradiction suppose the above statement is false. Then there exists $\varepsilon _0>0$, a sequence $\rho _n\to 0$ and $\psi _n\in \bigwedge ^2 (B_1^m)$ with
$$
\spt \psi_n \subset E_0^c\ ,
$$
where $E_0:=\{ x\in B_1^m\ :\ \dist {\cup _i D_i}\leq \varepsilon _0 |x|\}$, and yet
$$
C_{\rho _n, 0}(\psi _n)\neq 0\ .
$$
Then there exist $x_n\in E_0^c$ with $\lim _{r\to 0} M(C_{r, x_n})\neq 0$ and from the monotonicity formula we get that $M\big(C_{\rho _n |x_n|, 0}\res B_{\frac{\varepsilon _0}{2}}\big(\frac{x_n}{|x_n|}\big) \big)\geq \frac{\varepsilon _0^2}{4}\pi$. Taking a subsequence so that $\frac{x_n}{|x_n|}\to x_\infty$ (with $x_\infty \in E_0^c$) we have that 
$$
M\big(C_{\rho _n |x_n|, 0}\res B_{\frac{3\varepsilon _0}{4}}(x_\infty ) \big)\geq \frac{\varepsilon _0^2}{4}\pi\ .
$$
Since we also have
\bea
\nonumber \Big| C_{\rho _n |x_n|, 0}\res B_{\frac{3\varepsilon _0}{4}}(x_\infty ) \big( (\rho _n |x_n|)^2 \lambda ^{\rho _n |x_n|, 0^*}\omega -\omega _0 \big) \big) \Big| & & \\
\nonumber & \hspace{-12 cm} \leq & \hspace{-6 cm} M(C_{\rho _n |x_n|, 0}\res B_{\frac{3\varepsilon _0}{4}}(x_\infty)) \Big\| (\rho _n |x_n|)^2 \lambda ^{\rho _n |x_n|, 0^*}\omega -\omega _0 \Big\| _\infty \to 0\ ,
\eea
we conclude that we would get
$$
C_{\infty , 0}\res B_{\frac{3\varepsilon _0}{4}}(x_\infty)(\omega _0)\geq \frac{\varepsilon _0^2}{4}\pi\ ,
$$
which contradicts the fact that $C_{\infty, 0}$ is the union of $Q$ $\omega _0$-calibrated disks, i.e. has support in $E_0$. Hence the lemma holds true.


\bibliography{paper.bbl}
\bibliographystyle{plain}

\end{document}